\documentclass[11pt]{article}

\usepackage{amsmath}
\usepackage{amsthm}
\usepackage{amssymb}
\usepackage{amsfonts}

\newtheorem{Theorem}{Theorem}

\newtheorem{Lemma}[Theorem]{Lemma}
\newtheorem{Proposition}[Theorem]{Proposition}
\newtheorem{Corollary}[Theorem]{Corollary}

\newcommand{\R}{\mathbb{R}}
\newcommand{\Z}{\mathbb{Z}}
\newcommand{\E}{\mathbb{E}}
\newcommand{\I}{\mathbb{I}}
\newcommand{\Pp}{\mathbb{P}}
\newcommand{\N}{\mathbb{N}}
\newcommand{\Law}{\mathbb{L}}

\newcommand{\La}{\langle}
\newcommand{\Ra}{\rangle}

\newcommand{\llangle}{\langle\!\langle}
\newcommand{\rrangle}{\rangle\!\rangle}

\begin{document}

\title{$A+A \to A$, $\; \; B+A \to A$.}
 
\author{Roger Tribe
and Oleg Zaboronski
\\
Mathematics Institute, University of Warwick,
Coventry, CV4 7AL, UK}

\maketitle

\begin{abstract}
This paper considers the decay in particle intensities for a translation invariant two species system of diffusing
and reacting particles on $\Z^d$ for $d \geq 3$. The intensities are shown to
approximately solve modified rate equations, from which their polynomial decay 
can be deduced. The system illustrates that the underlying diffusion and reaction rates
can influence the exact polynomial decay rates, despite the system evolving in
a supercritical dimension. 
\end{abstract}
\section{Introduction} \label{s1}
Our aim is to analyse the decay of densities in the two species reacting system $A+A \to A$
and $A+B \to A$, when particles perform simple random walks on $\Z^d$ for $d \geq 3$. 
At time  zero particles are, for example, Poisson with uniform (non-zero) intensities.  
The density of $A$ particles should fall as $O(t^{-1})$, for typical random walks and reaction mechanisms. Our aim here is to show that the density of $B$ particles falls like $O(t^{-\theta})$, for some $\theta$ to be found and which depends on the details of the random walks and reactions. 

We specify one specific set of reactions. Type $A$ particles perform rate $D_A$ independent simple random walks on $\Z^d$, for $d \geq 3$. In addition, 
each pair of type $A$ particles on the same site will (independently of other pairs) coalesce at rate $2 \lambda_A$.  
Type $B$ particles perform rate $D_B$ simple random walks on $\Z^d$, independently of the $A$ population. Moreover 
any pair of one $A$ particle and one $B$ particle at the same site, will (independently of all other randomness) coalesce into a type $A$ particle at rate $\lambda_B$.  
Note that the $A$ population is unaffected by the presence of the $B$ population.
All rates are strictly positive.

Suppose at $t=0$ the numbers of $A$ and $B$ particles are independent and identically distributed
at sites $x \in \Z^d$, with finite moments of all orders. 
The methods of Bramson and Griffiths \cite{bramson1980asymptotics} and of van den Berg and Kesten \cite{kesten1998asymptotic,van2002randomly} 
imply the asymptotics 
\begin{equation} \label{asymptoticA}
\Pp[ \mbox{there is an $A$ particle at the origin}] \sim \frac{1}{p_A \lambda_A t}
\quad \mbox{as $t \to \infty$}
\end{equation}
where 
\[
p_A = \frac{\gamma D_A }{\gamma D_A + \lambda_A}
\] 
is the probability that a system with just two $A$ particles, both starting at the origin, never coalesces, 
and $\gamma$ the probability that a simple random walk starting at $0$ never makes a return visit to the origin.
Note the constant $\frac{1}{p_A \lambda_A}$ holds universally over the possible initial conditions. 

The aim of this paper is to show that 
\[
\Pp[ \mbox{there is an $B$ particle at the origin}] \sim  C \; t^{- \frac{p_B \lambda_B}{p_A \lambda_A}}
\quad \mbox{as $t \to \infty$}
\]
where $C=C(D_A,D_B,\lambda_A,\lambda_B)>0$ and 
\[
p_B= \frac{\gamma (D_A+D_B) }{\gamma (D_A+D_B) + \lambda_B}
\]
is the probability that a system with just one type $A$ and one type $B$ particle, both starting at the origin, 
never coalesces. This leads to the decay rate
$C t^{-\theta}$ for the $B$ population where
\[  
\theta = \frac{D_A+D_B}{D_A} 
\frac{1+\gamma D_A \lambda_A^{-1}}{1+\gamma(D_A+D_B)\lambda_B^{-1}}.
\]
Note that this $\theta$ can take any value in $(0,\infty)$ and also the good sense that
$\theta \to 0$ as $\lambda_B \to 0$ and $\theta \to \infty$ as $\lambda_A \to 0$.

Taking $\lambda_A \to \infty$ or $\lambda_B \to \infty$ we expect to get the decay rates for
the case of instantaneous coalescence. In particular 
\begin{equation} \label{instantly}
p_A \lambda_A \to  \gamma D_A, \quad p_B \lambda_B \to \gamma(D_A+D_B) 
 \quad \mbox{as $\lambda_A, \lambda_B \to \infty$}
\end{equation}
so that $\theta \to 1+ (D_B/D_A)$. Notice that this limiting answer for $\theta$  has a higher degree universality than the exponent for finite reaction rates; the limit depends on the ratio $D_B/D_A$, but not the constant $\gamma$, which depends on the geometry of the underlying lattice. Therefore we conjecture that the infinite reaction rate formula $\theta=1+ (D_B/D_A)$ is valid for a large class of translation invariant lattices (for example the triangular lattice).

These results should be compared with the solutions to naive rate equations:
\begin{equation} \label{naiveRE}
\frac{da}{dt} = - \lambda_A a^2, \quad \frac{db}{dt} = - \lambda_B ab,
\end{equation}
which yield  $a(t) \sim 1/\lambda_A t$ and $b(t) \sim C
\; t^{-\lambda_B/\lambda_A}$ as $t \to \infty$. Thus the naive rate 
equations correctly predict the polynomial rate of decay for $A$ particles
(with an incorrect constant) and fail to predict the correct decay rate for $B$ particles.
 The point is that, in supercritical dimensions,  
such rate equations do give the right answers providing the
constants are modified to 
\begin{equation} \label{modifiedRE}
 \frac{da}{dt} = - p_A \lambda_A a^2, \quad \frac{db}{dt} = - p_B 
\lambda_B ab.
\end{equation}
Since the constants for this system of two equations get 
involved in the polynomial decay rate for $B$ particles, one cannot read off the correct decay rates
just from the naive rate equations and the work establishing the modified constants is 
crucial. 
Our system was chosen as the simplest system to illustrate this point.
This method, known as the 'modified rate equation method',
should apply
in all transient dimensions $d \geq 3$. 

Our approach follows that in Van den Berg and Kesten \cite{kesten1998asymptotic,van2002randomly} on single species coalescent models, where enough asymptotic independence of sites 
is established to show that the modified rate equations give the correct decay rates. However, unlike \cite{kesten1998asymptotic,van2002randomly},
we will exploit negative dependence inequalities for the $A$ particle density. Van den Berg and Kesten show that
these can be avoided, estimating covariances directly, and 
hence their results apply to a class of related coalescence models, with more complicated coalescence mechanisms, for which negative correlation is 
unknown.  In our two species setting, the $B$ particle decay rate can be arbitrarily slow which necessitates more work 
in the error estimates. We use high moments to establish 
the de-correlation of the $A$ and $B$ particle densities, which in turn exploit the 
negative dependence properties of the $A$ particles. To estimate these high moments, we use 
a Marcinkiewicz-Zygmund type inequality (also called a square function inequality)
for variables with negative dependence, see Proposition \ref{MZ}, which we believe is new and which may be of independent interest. 

\noindent
\textbf{Related literature.} 
The two-species reaction $A+A \to A, \; A+B \to A$ 
has been studied in theoretical physics in the context of
persistence, where $B$'s play the role of test particles in the environment of reacting $A$ particles. Then the density of $B$ particles at time $t$ reflects the probability that a test particle has not collided with an $A$ particle through time $t$.  In particular, in \cite{redner} the model
has been treated in all spatial dimensions using a combination of rate and Smoluchowski equations;
in \cite{monthus} the one-dimensional case was analyzed by developing a perturbation theory
in the small mobility of $B$ particles; in \cite{howard}, \cite{rajesh} the perturbative renormalisation group
method was applied to calculate the answers for the space dimension $d=2-\epsilon$, where $\epsilon\geq 0$. 
(The final answers presented in \cite{howard} are incorrect even within the non-rigorous renormalisation
group framework, as their derivations neglect the anomalous dimension of the field operator describing the
$B$ particles.)

We briefly discuss work in $d=1,2$. 
In the critical dimension $d=2$, the modified rate equation method was shown 
to work for a single species coalescent system in \cite{lukins2018multi}, to derive the 
leading asymptotic of the decay rates of all multi-point densities,
where logarithms modifications appear to the $d \geq 3$ decay rates, 
and the leading asymptotic $a(t) \sim C \log(t)/t $ 
should be universal over a range of reaction mechanisms due to the recurrence.
The two species system in $d=2$ was studied using non-rigourous renormalisation group 
calculations in \cite{rajesh} and this suggests a surprising answer for the $B$ particle 
density decay, where the logarithmic modifications to the answers are non-universal, 
in the sense that the powers of $\log t$ appearing in a full asymptotic expression for $b(t)$ 
depend on both the reaction and diffusion rates. 
For coalescing systems in dimension $d=1$ the method fails, and the multi-point densities are more complicated. Luckily, for instantaneous coalescence,  techniques from integrable probability allow all multi-point densities to be calculated, see \cite{tz3}, \cite{garrod2018examples}.

\noindent
\textbf{Layout.} $\,$
In section \ref{s2} we give a construction of the particle systems, explain the main intuition of the proof, 
following the arguments from  \cite{kesten1998asymptotic,van2002randomly}, 
and state the main result. 
In section \ref{s3} we establish a negative dependence property of the $A$ particle density, 
and prove the Marcinkiewicz-Zygmund type inequality.
We complete the proof of the main result in section \ref{s4}. 

\noindent
\textbf{Notation.} $\,$
We will use small case constants $c_1(p), c_2(t) \ldots$ for constants that will be referred 
to later in the paper. We use large case constants $C_{p,t,\ldots}$ for constants
whose exact value is unimportant and may vary from line to line. In all cases the dependence
on all variables will be indicated, except that we suppress dependence on the dimension $d$.  
\section{Construction} \label{s2}
\subsection{A strong equation}
We now give a  construction for the processes with finite reaction rates $\lambda_A, \lambda_B < \infty$.
We construct the process as the solution of a system of equations on $\Z^d$, driven by
Poisson processes (very close to standard graphical constructions of particle systems). 
We start with a probability space equipped with a collection of Poisson families as follows.
Write $x \sim y$, for $x, y \in \Z^d$, to mean that they are nearest neighbours, that is that
$|x-y|=1$ (in the Euclidean norm).
\begin{itemize}
\item 
The family $(P^A(i,x,y): i \in \N, \, x,y \in \Z^d, \, x \sim y)$ of i.i.d. rate $D_A/2d$ Poisson processes 
will control the jumps of the $i$th type $A$ particle at $x$ to the site $y$. 
\item
The family $(P^B(i,x,y): i \in \N, \, x,y \in \Z^d, \, x \sim y)$ of i.i.d. rate $D_B/2d$ Poisson processes 
will control the jumps of the $i$th type $B$ particle at $x$ to the site $y$. 
\item
The i.i.d. Poisson family  $(P^{AA}(i,j,x): i,j \in \N, \, x \in \Z^d)$, at rate $\lambda_A$, will control the 
coalescence of the $i$th type $A$ particle onto the $j$th type $A$ particle at site $x$ (the
total rate of coalesence for each pair will be $2 \lambda_A$).  
\item
The i.i.d. Poisson family  $(P^{AB}(i,j,x): i,j \in \N, \, x \in \Z^d)$, at rate $\lambda_B$, will control the 
coalescence of the $i$th type $A$ particle with the $j$th type $B$ particle at site $x$. 
\end{itemize}
These families are  independent of each other, are defined on a filtered probability space
$(\Omega, \mathcal{F}, (\mathcal{F}_t), \Pp)$, and are compatible with the filtration, that is that they are adapted and 
increments over any time interval $[s,t]$ are $\mathcal{F}_s$ independent.

The processes $\xi= (\xi_t(x): t \geq 0, x \in \Z^d)$ and $\eta= (\eta_t(x): t \geq 0, x \in Z^d)$ will record the 
particle numbers at sites $x \in \Z^d$ and times $t \geq 0$ for the $A$ and $B$ type particles. 
They solve, for all $x \in \Z^d$, the stochastic differential equations
\begin{eqnarray}
d \xi_t(x) &=&   \sum_{y:y \sim x} \sum_{i \geq 1} \left( 
 I(i \leq \xi_{t-}(y)) dP^{A}_t(i,y,x) - I(i \leq \xi_{t-}(x)) dP^A_t(i,x,y)\right)
 \nonumber \\
&& \hspace{.7in}  -  \sum_{i,j \geq 1} I(i \vee j \leq \xi_{t-}(x), i \neq j) dP^{AA}_t(i,j,x), \quad
\mbox{for $t \geq 0$, \, a.s.,} \label{xieqn} \\
d \eta_t(x) &=&  \sum_{y:y \sim x} \sum_{i \geq 1}\left( 
I(i \leq \eta_{t-}(y)) dP^B_t(i,y,x) - I(i \leq \eta_{t-}(x)) dP^B_t(i,x,y) \right)
\nonumber \\
&& \hspace{.4in}  -  \sum_{i,j \geq 1} I(i \leq \xi_{t-}(x), j \leq \eta_{t-}(x)) dP^{AB}_t(i,j,x), \quad
\mbox{for $t \geq 0$, \, a.s.} \label{etaeqn}
\end{eqnarray}

To formulate an existence and uniqueness result, we will look for solutions that live in certain 
weighted $l_{2}(\Z^d)$ spaces. 
We will be interested in translation invariant solutions and the exact choice of weights 
is not so important; to be concrete
we define the space of functions with slower than exponential growth by
\[
S_{\alpha} = \{f: \Z^d \to \N \, | \, \|f\|_{\alpha} := \sum_x |f(x)|^2 e^{- \alpha |x|} < \infty\} \quad \mbox{for $\alpha>0$.}
\]
With the norm $\|f\|_{\alpha}$, the space  $S_{\alpha}$ is a complete metrisable space, $S_{\alpha}^2$ is the product space and $D_{S_{\alpha}^2}[0,\infty)$ the space of 
cadlag $S_{\alpha}^2$ valued paths in the Skorokhod topology.

%
%
\begin{Proposition} \label{strong solutions}
Suppose $\xi_0$ and $\eta_0$ are $\mathcal{F}_0$ measurable variables satisfying
$\E[\|\xi_0\|_{\alpha}] + \E[\|\eta_0\|_{\alpha}] <\infty$ for some $\alpha>0$. 
Then there exist an adapted solutions $(\xi_t,\eta_t)_{t \geq 0}$ to (\ref{xieqn}) and (\ref{etaeqn}) with 
paths in $D_{S^2_{\alpha}}[0,\infty)$. 
Solutions are pathwise unique and unique in law. If
$(\xi_0,\eta_0)$ are translation invariant then so too are $(\xi_t,\eta_t)$ at all $t \geq 0$. 
Whenever $\E[ \sum_x (|\xi_0(x)|^p + |\eta_0(x)|^p) e^{- \alpha |x|}] < \infty$, for some $p>0$, then 
the same moments remain bounded over finite time intervals.
\end{Proposition}

\noindent
The result can be established by, for example, a standard Picard scheme and 
Gr\"{o}nwall estimates for the moments. Similar infinite systems of Poisson driven equations are discussed in, for example, 
\cite{graham1996weak} in the context of spin systems. 
\subsection{The main result} \label{s2.2}
\textbf{Initial Conditions.}
For initial conditions we ask that $(\xi_0(x): x \in \Z^d)$ and $(\eta_0(x):x \in \Z^d)$ are independent i.i.d. families, and that $\E[ \xi_0^p(0) + \eta^p_0(0)] < \infty$
for all $p \geq 0$. 
We will also place some restriction on the law of $\xi_0(0)$ 
so that we can use BKR inequalities to derive negative dependence properties. 
We ask that the law of $\xi_0(0)$ lies in 
\begin{equation} \label{initialcondition}
\overline{\mathcal{M}_{PB}} = \mbox{closure of the set of Poisson-Binomial variables}
\end{equation}
where the closure is in the topology of convergence in distribution, and a Poisson-Binomial
variables is a variable $X$ that can be written as a sum $X = \sum_{i=1}^N B(p_i)$ of 
independent Bernouilli $B(p_i)$ variables.  In particular deterministic $\xi_0$
or Poisson $\xi_0$ are special cases. 

\begin{Theorem} \label{mainresult}
Suppose that $(\xi,\eta)$ are solutions to (\ref{xieqn}) and (\ref{etaeqn}) in $d \geq 3$
with non-zero initial conditions as above.
Then there exists $\kappa_1= \kappa_1(d) >0$ so that
\[
\E[ \xi_t(0)]  =  \frac{1}{p_A \lambda_A} t^{-1} \left(1+ O(t^{-\kappa_1})\right)
\]
and there exists $ \kappa_2= \kappa_2(d) >0 $ and $c_0$ (where $c_0$ may depend on
the initial conditions and all rates $D_A,D_B,\lambda_A,\lambda_B$)
so that
\[
\E[ \eta_t(0)] = c_0 t^{-\frac{p_B \lambda_B}{p_A \lambda_A}} \left(1+ O(t^{-\kappa_2})\right).
\]
\end{Theorem}

\noindent
\textbf{Remark 2.2.1} 
In the instantly coalescing cases,
where either $\lambda_A$ or $\lambda_B$ or both are infinite,
we believe that the same results holds except that either or both of the 
substitutions in (\ref{instantly}) are required. 
We discuss two methods to establish this belief at the end of Section \ref{s4}. 

\noindent
\textbf{Remark 2.2.2}
The leading asymptotics can be stated in terms of occupation probabilities densities since 
\[
\Pp[ \xi_t(0)>0] \sim \E[ \xi_t(0)] \quad \mbox{and} \quad \Pp[ \eta_t(0)>0] \sim \E[ \eta_t(0)] \quad \mbox{as $t \to \infty$.}
\]

\noindent
\textbf{Remark 2.2.2} Concrete values for $\kappa_1,\kappa_2$, which we do not claim are optimal,
can be read off from the proofs; for example any $\kappa_1 < (d-2)/(3d-4)$ is achievable (see Remark \textbf{4.2.1}); one can take
$\kappa_2 = \frac{1}{17}$ in $d=3$  and $\kappa_2 = \kappa_1$ for large $d$ (see Remark \textbf{4.4.1}).

In the rest of this subsection we repeat the heuristic argument from \cite{kesten1998asymptotic} that guide 
the calculations, adapting it for our two species case. We can divide the proof into three steps.

\noindent
\textit{Step I.} An exact calculation shows that
\[
\frac{d}{dt} \E[\eta_t(0)] = - \lambda_B \E[\xi_t(0)\eta(0)].
\]
At large times, the expectation $\E[\xi_t(0)\eta(0)]$ is well approximated by $\Pp[ \xi_t(0) = \eta_t(0) =1]$
since the occurrence of sites with three or more 
particles will be of smaller order in $t^{-1}$. 
Trace an A and a B particle that are at $0$ at time $t$ back to the time $t-s$, where
we will take $1 \ll s \ll t$. These particles must have been located at some points $x$ and $y$, 
suggesting that (using an informal notation $\approx$ for 'approximately equal to')
\[
\Pp[ \xi_t(0) = \eta_t(0) = 1] 
 \approx   \sum_{x,y} \Pp[ \xi_{t-s}(x)=1, \eta_{t-s}(y)=1] \psi_s(x,y) 
\]
where $\psi_s(x,y)$ is the probability an A and a B particle starting at $x$ and $y$ end up at $0$
at time $s$ without coalescing. Step I is to justify (and quantify) these approximations, which all 
follow from the density of particles being low, requiring $t-s \gg 0$.

\noindent
\textit{Step II.}  A random walk calculation shows that, 
when $1 \ll |x-y|= O(t^{1/2})$ and $s \gg 1$, 
\[
\psi_s(x,y) \approx p_B \, p^A_{s}(x) p^B_{s}(y) 
\]
where $p^A_s(x)$ and $p^B_s(x)$ are the simple random walk transition densities for
$A$ and $B$ particles. This is intuitive if one reverses the paths and look at them 
starting at $0$ and reaching $x$ and $y$ without coalescing, since then coalescence in 
transient dimensions is mostly decided by times smaller than $s$. Combining this with step I leads to 
\[
\frac{d}{dt} \E[\eta_t(0)]
\approx - \lambda_B p_B \E \left[ \sum_{x,y} \xi_{t-s} (x) p^A_s(x) \eta_{t-s}(y) p^B_s(y) \right].
\]
The estimates for this approximation step will be taken from \cite{kesten1998asymptotic} Lemma 12.
 
 \noindent
\textit{Step III.} The hardest step is to get some quantitative control on the approximate independence
\[
\Pp[ \xi_{t-s}(x)=1, \eta_{t-s}(y)=1] \approx \Pp[ \xi_{t-s}(x)=1] \Pp[ \xi_{t-s}(y)=1]
\]
which holds at least for $|x-y| \gg 1$ and $t-s \gg 1$. We bound the covariance of $\sum \xi_t(x) f(x)$ and 
$ \sum \eta_t(x) g(x)$ using H\"{o}lder's inequality and high fluctuation moments for the $A$ population 
(which will exploit the negative dependence properties). This will 
be sufficiently to justify the final step 
\begin{eqnarray*}
\frac{d}{dt} \E[\eta_t(0)] & \approx &  - \lambda_B p_B 
\E \left[ \sum_x \xi_{t-s} (x) p^A_s(x) \right] \E \left[ \sum_x \eta_{t-s}(x) p^B_s(x) \right] \\
& = &  - \lambda_B p_B  \E [ \xi_{t-s} (0)] \E [\eta_{t-s}(0)] \\
&  \approx &  - \lambda_B p_B \E [ \xi_{t} (0)] \E [\eta_{t}(0)].
\end{eqnarray*}
This leads to the second of the modified rate equations (\ref{modifiedRE}) (and a similar and simpler 
argument will justify the first of the modified rate equations). The resulting 
asymptotics follow from analysing these equations and the error estimates in the approximations above.

\section{Negative dependence} \label{s3}
Negative dependence properties for the one species coalescing system were first used in
Arratia \cite{arratia1981limiting} to estimate variances, and subsequently in 
\cite{kesten1998asymptotic,van2002randomly},  \cite{lukins2018multi}.
The property of negative association, and the weaker property 
LNQD, are some of a number of related definitions - see Newman \cite{newman1984asymptotic} or Pemantle \cite{pemantle2000towards} for discussions.

We will broadly follow the van den Berg and 
Kesten argument from \cite{van2002randomly}, where they deduce a negative dependence property 
for some coalescing systems from the van den Berg-Kesten-Reimers (BKR)
inequality. In particular, we modify the nice technique used in \cite{van2002randomly} of obtaining our desired system 
as the limit of projections of suitable multi-colour particle systems. 

In section \ref{s3.1} we state the 
main negative correlation property for the $A$-particles, and derive consequences for moments of the process. 
In section \ref{s3.2} we derive a 'square function' inequality (analogous to half of a 
Burkh\"{o}lder inequality) for variables with a suitable negative dependency. 
The proof of negative association, following closely that for the systems in
\cite{van2002randomly}, is given in \ref{s3.3}.
%
%
%
\subsection{Negative association for $A$ particles} \label{s3.1}
%
%
%
\noindent
\textbf{Definition} 
A random vector $X=(X_1,\ldots,X_N) \in \R^N$ is said to be {\it negatively associated} if
\[
\E \left[ f(X_{i_1},\ldots,X_{i_k}) g(X_{j_1},\ldots,X_{j_l}) \right] \leq 
\E \left[ f(X_{i_1},\ldots,X_{i_k})] \E[ g(X_{j_1},\ldots,X_{j_l}) \right] 
\]
whenever $\{i_1,\ldots,i_k\}$ and $\{j_1, \ldots,j_l\}$ are disjoint subsets of $\{1,\ldots,N\}$
and $f:\R^k \to \R$ and $g:\R^l \to \R$ are bounded increasing functions. 
An infinite vector is called negatively associated, if every finite subvector
has the corresponding property. 
\begin{Proposition} \label{NA}
Let $\xi$ be the solution of the $A$-particle equation (\ref{xieqn}) started from 
an initial condition $(\xi_0(x):x \in \Z^d)$ of independent variables from
$\overline{\mathcal{M}_{PB}}$. Then for any $t>0$ the vector
$(\xi_t(x): x \in \Z^d)$ is negatively associated. In addition 
it satisfies the inequality
\begin{equation} \label{extradependence}
\Pp \left[ \xi_t(x) \geq k+l \right] \leq \Pp \left[ \xi_t(x) \geq k \right] \, \Pp \left[ \xi_t(x) \geq l \right]
\quad \mbox{for $k,l \in \N$ and $x \in \Z^d$.}
\end{equation}
\end{Proposition}
The proof is given in section \ref{s3.3}.
\begin{Corollary} \label{moments1}
For $\xi$ as in Proposition \ref{NA} we have
\begin{equation} \label{naone}
\E \left[ \xi_t(x_1) \ldots \xi_t(x_n) \right] \leq \E \left[ \xi_t(x_1)\right] \ \ldots  \E \left[  \xi_t(x_n) \right]
\; \mbox{for disjoint $x_1,\ldots,x_n$ and
$t \geq 0$,}
\end{equation}
and
\begin{equation} \label{natwo}
\E \left[ \xi_t(x) (\xi_t(x)-1) \ldots (\xi_t(x)-(n-1))\right] \leq n! \left( \E \left[ \xi_t(x) \right]\right)^n 
\quad \mbox{for $x \in \Z^d, \, n \in N$ and $t>0$.}
\end{equation}
\end{Corollary}
\noindent
\textbf{Proof of Corollary \ref{moments1}.} 
By the negative association and induction in $n$ we have 
\[
\E \left[ \prod_{i \leq n} \left(\xi_t(x_i) \wedge K\right) \right] \leq \prod_{i \leq n} \E \left[ \xi_t(x_i) \wedge K \right] 
\]
and letting $K \uparrow \infty$ yields (\ref{naone}). Starting with (\ref{extradependence}), induction 
yields, for $k_1, \ldots, k_n \in \N$,
\[
\Pp \left[ \xi_t(x) \geq k_1 + \ldots + k_n \right] \leq \Pp \left[ \xi_t(x) \geq k_1 \right] \ldots 
\Pp \left[ \xi_t(x) \geq k_n  \right]
\]
Now (\ref{natwo}) follows from the identity
\[
X (X-1) \ldots (X-(n-1)) = n! \sum_{i_1=1}^{\infty} \ldots 
 \sum_{i_n=1}^{\infty} \I(X \geq i_1+ \ldots + i_n).
\]
\qed

\noindent
\textbf{Remark 3.1.1}
The properties (\ref{extradependence}) and (\ref{natwo}) are instances of the factorial moment negative dependence 
discussed in Liggett \cite{key1957219m} for systems of independent random walks. However the method of proof used in \cite{key1957219m} does not seem to extend to establish such properties for coalescing systems.

\noindent
\textbf{Remark 3.1.2} We have not been able to establish useful negative dependence properties between A and B particles.
As a guiding example, consider the two processes started with exactly one A particle and one B particle, both at the origin. 
Then $\E[\xi_t(x) \eta_t(0)] \leq \E[\xi_t(x)] \E[\eta_t(0)] $ for $|x|$ small, but
the reverse inequality holds for $|x|$ large.

\subsection{A Marcinkiewicz-Zygmund inequality for negatively dependent
variables} \label{s3.2}
The following is an analogue of the classical Marcinkiewicz-Zygmund inequality
for independent variables, or Burkh\"{o}lder type inequality for martingales
(see \cite{chow2003probability}). It only uses a slightly weaker form of 
negative correlation called LNQD (linearly negative quadrant dependant - see \cite{newman1984asymptotic}) 
which we now recall. 

\noindent
\textbf{Definition} 
A random vector $X=(X_1,\ldots,X_N) \in \R^N$ is called {\it linearly negative quadrant dependent} (LNQD) if the pair 
$\sum_{i \in A} \lambda_i X_i$ and $\sum_{j \in B} \mu_j X_j$ are negatively associated 
for any disjoint $A$ and $B$ and positive $\lambda_i, \mu_j$. 
An infinite vector is called LNQD if every finite subvector has the corresponding property. 

\begin{Proposition} \label{MZ}
For $p \in \N$ there exists constants $C_p< \infty$ so that 
if $X=(X_1, \ldots,X_N)$ is LNQD for some $N \geq 1$, $E[X_i]=0$ and $E[|X_i|^{2p}] < \infty$ for $i=1,\ldots,N$, then
\[
\E \left[ \left| \sum_{i=1}^N X_i \right|^{2p} \right] \leq C_p \, \E \left[  \left| \sum_{i=1}^N X_i^2 \right|^{p} \right].
\]
\end{Proposition}

\textbf{Proof.} Unindexed sums $\sum_i$ will be over the indices $i =1,\ldots,N$. 
Set $S= \sum_i X_i$ and $S_i = \sum_{j>i}^N X_j$ and $[S] = \sum_i X^2_i$. 
The binomial expansion 
\[
S^{2p} = (X_1+S_1)^{2p} = \sum_{k=0}^{2p} {2p \choose k} X_1^{k} S_1^{2p-k} =
S^{2p}_1 + \sum_{k=1}^{2p} {2p \choose k} X_1^{k} S_1^{2p-k} 
\]
may be iterated, be repeating the expansion on the terms $S_1^{2p}, S_2^{2p}, \ldots, S_{N-1}^{2p}$, to find
\begin{equation} \label{expansion}
S^{2p} = S_{N-1}^{2p} + \sum_{k=1}^{2p} {2p \choose k} \sum_{i=1}^{N-1} X_i^{k} S_i^{2p-k}
=  \sum_{k=1}^{2p} {2p \choose k} \sum_i X_i^{k} S_i^{2p-k}
\end{equation}
where in the last equality we are using the convention that $S_N^k = 0$ for $k \geq 1$ but $S_N^0=1$. 
Taking expectations, the term in (\ref{expansion}) when $k=1$, that is $\sum_i X_i S_i^{2p-1}$,
has negative expectation by the LNQD assumption; indeed we may apply the negative association to the 
increasing functions $f(x)=x$ and $g(x) = x^{2p-1}$ to the truncated variables $\psi_K(X_i)$ and $\psi_K(S_i)$
for $\phi_K(x) = (x \wedge K) \vee (-K)$ and then let $K \uparrow \infty$ using the moment assumption to remove
the truncation.  
 Also the term in (\ref{expansion}) when $k=2p$ satisfies $\sum_i X_i^{2p} \leq (\sum_i X_i^2)^p = [S]^p$. 
The aim is to control the remaining terms in (\ref{expansion}) using $[S]$.

Consider the term $\sum_i X_i^k S_i^{2p-k}$ for some $k \in \{2,3,\ldots,2p-1\}$. 
In the martingale proof (that is for a Burkh\"{o}lder inequality) one could first 
bound this by H\"{o}lder's inequality, for any $\delta >0$, by 
\begin{eqnarray}
\E \left[ \sum_i X_i^k S_i^{2p-k} \right] & \leq &  \E \left[ \sup_{j \leq N} |S_j|^{2p-k} \sum_i |X_i|^k \right]
\nonumber \\
& \leq &  \frac{2p-k}{2p} \delta \, \E \left[ \sup_{j \leq N} |S_j|^{2p} \right]  +
\frac{k}{2p} \delta^{-(2p-k)/k} \, \E \left[ ( \sum_i |X_i|^k )^{2p/k} \right] \nonumber \\
& \leq &  C_{p,k} \delta \, \E[ \sup_{j \leq N} |S_j|^{2p}]  +
C_{p,k,\delta} \, \E[ [S]^p]. \label{mgcase}
\end{eqnarray}
Then Doob's maximal inequality allows one to bound $ \E[ \sup_j |S_j|^{2p}]$ by $ C_p \E[S^{2p}]$. 
Choosing small $\delta$ one finds
\[
\E[S^{2p}] \leq C_p \E[[S]^p] + (1/2) \E[S^{2p}]
\] 
and one can close the inequality. In our case we don't see how to get a suitable 
maximal inequality working just from a negative dependence assumption. 

However there is alternative argument. If we symmetrize (\ref{expansion}) 
over all $N!$ possible orderings of the variables $X_1,\ldots,X_N$ we obtain 
\begin{equation} \label{1*}
S^{2p} = \sum_{k=1}^{2p} {2p \choose k} \sum_i X_i^{k} Z(i,k)
\end{equation}
where $Z(i,2p)=1$ and for $k =1,\ldots,2p-1$
\[
Z(i,k) = \frac{1}{N} \sum_{j=1}^{N-1} {N-1 \choose j}^{-1}
\sum_{\Lambda \subseteq \{1,\ldots,N\} \setminus \{i\}, |\Lambda|=j } S^{2p-k}_{\Lambda}
\]
with $S_{\Lambda}= \sum_{i \in \Lambda} X_i$. 
 Again, the term $k=1$ has negative expectation by  
the LNQD assumption, so we need we need to 
bound the terms when $k=2,\ldots,2p-1$. 
Let $A$ be a random subset of $\{1,\ldots,N\}$ chosen as follows. First choose 
$K$ uniformly from $\{0,1,\ldots,N\}$ and then, conditional on $K$, choose $A$ uniformly from all
subsets of $\{1,\ldots,N\}$ of size $K$ (so that $A = \emptyset$ if $K=0$).
Then, for a fixed $i$ and fixed $\Lambda  \subseteq \{1,\ldots,N\} \setminus \{i\}$ with $|\Lambda|=j$
\begin{eqnarray*}
\Pp[A \setminus \{i\} = \Lambda] &=&  \Pp[A= \Lambda] + \Pp[A= \Lambda \cup \{i\}] \\
& = & \frac{1}{N+1} \left( {N \choose j}^{-1} + {N \choose j+1}^{-1} \right) \\
& = & \frac{1}{N} {N-1 \choose j}^{-1}.
\end{eqnarray*}
We may assume (by extending the probability space if necessary) that $(K,A)$ are defined on the same probability
space as $X$ but are independent of $X$. 
Then the term $\E[\sum_i X_i^k Z(i,k)]$ can be rewritten (with the convention that $S_{\emptyset}=0$) as 
\begin{equation} \label{2*}
\E \left[ \sum_i X^k_i Z(i,k) \right] = \E \left[ \sum_i X_i^{k} S_{A \setminus \{i\}}^{2p-k} \right].
\end{equation}
The idea is to approximate each term $S_{A \setminus \{i\}}$ by $S_A$, a common factor independent of $i$.
Since
\[
|S_{A \setminus \{i\}}|^{2p-k} \leq C_{p,k} \left( |S_A|^{2p-k} + |X_i|^{2p-k} \right)
\]
we have the bound, arguing as in (\ref{mgcase}), for any $ k \in \{2,3,\ldots,p-1\}$ and $\delta>0$
\begin{eqnarray}
\E \left[ \sum_i X_i^{k} S_{A \setminus \{i\}}^{2p-k} \right] & \leq &
C_{p,k} \E \left[ |S_{A}|^{2p-k} \sum_i |X_i|^{k} \right] + C_{p,k} 
\E \left[ \sum_i |X_i|^{2p} \right] \nonumber \\
& \leq & C_{p,k} \, \delta \, \E[ S_A^{2p}]  + C_{p,k,\delta} \, \E[[S]^{p}]. \label{3*}
\end{eqnarray}
To close the inequality it remains to show that we can bound $\E[S_A^{2p}]$ in 
terms of $\E[S^{2p}]$. 
For large $N$, $S_A^{2p}$ should be well approximated by $ | \sum_i W_i X_i |^{2p}$   
where the vector $W=(W_1, \ldots,W_N)$ is
chosen independently of $X$ by first picking $K$ uniformly as above and then, conditionally on $K$, 
letting $W_i$ be i.i.d. Bernouilli $(K/N)$ variables. We  
claim (and prove below) that
\begin{equation} \label{claim10}
\Pp \left[ \sum_i W_i = j \right] \geq \frac{1}{4(N+1)} \quad \mbox{for $j=0,\ldots,N$.}
\end{equation}
Conditional on $\sum W_i = j$ the vector $(W_1, \ldots,W_N)$ is uniformly distributed over all vectors
with $j$ ones $N-j$ zeros. This and the claim (\ref{claim10}) now imply that
\begin{equation} \label{A*}
\E[ S_A^{2p}] \leq 4 \E \left[ | \sum_i W_i X_i |^{2p} \right]. 
\end{equation}
We can bound this easier sum by 
\begin{eqnarray}
\E \left[ | \sum_i W_i X_i |^{2p} \right]  
& \leq & 2^{2p-1} \E \left[ 
| \sum_i (W_i-(K/N)) X_i |^{2p} \right] + 2^{2p-1} \E \left[ (K/N)^{2p} | \sum_i X_i |^{2p} \right] \nonumber \\
& \leq & C_p \E \left[ 
| \sum_i (W_i-(K/N)) X_i |^{2p} \right] + C_p \E[S^{2p}].\label{B*}
\end{eqnarray}
But we can now use the case of Zygmund's inequality for independent variables to see
\begin{eqnarray}
 \E \left[ \left. | \sum_i (W_i-(K/N)) X_i |^{2p} \right| \; \sigma\{X,K\} \right]
& \leq & C_p  \E \left[ \left. | \sum_i (W_i-(K/N))^2 X_i^2 |^{p} \right| \; \sigma \{X,K\} \right] \nonumber \\
& \leq & C_p | \sum_i X_i^2 |^{p} = C_p [S]^p. \label{C*}
\end{eqnarray}
The last three steps (\ref{A*},\ref{B*},\ref{C*})  allow us to bound 
\[
 \E[ S_A^{2p}] \leq 4 \E[ | \sum_i W_i X_i |^{2p}] \leq C_p \E[S^{2p}] + C_p \E[[S]^p]
 \]
which along with (\ref{1*},\ref{2*},\ref{3*}) completes the proof by closing the inequality as in the martingale case. 

To establish the claim (\ref{claim10}) note that
\[
\Pp \left[ \sum_i W_i =n \right]  =  \frac{1}{N+1} {N \choose n} \sum_{l=0}^N \left(\frac{l}{N}\right)^n 
\left(1- \frac{l}{N}\right)^{N-n}.
\]
For $n=0$, or for $n=N$, this immediately gives the desired lower bound, by examining 
just the term $l=0$, or the term $l=N$.
When $n \in \{1,\ldots,N-1\}$ note that the function $x^n (1-x)^{N-n}$
is unimodal and achieves its maximum at $x=n/N$. 
Examining Riemann block diagrams one has
\begin{eqnarray*}
\int^{n/N}_0 x^n(1-x)^{N-n} dx & \leq & \frac{1}{N} \sum_{l=0}^n \left(\frac{l}{N}\right)^n \left(1- \frac{l}{N}\right)^{N-n}, \\
\int_{n/N}^1 x^n(1-x)^{N-n} dx & \leq & \frac{1}{N} \sum_{l=n}^{N} \left(\frac{l}{N}\right)^n 
\left(1- \frac{l}{N}\right)^{N-n}
\end{eqnarray*}
and so 
\[
\sum_{l=0}^N \left(\frac{l}{N}\right)^n \left(1- \frac{l}{N}\right)^{N-n} \geq (N/2) \int^1_0 x^n(1-x)^{N-n} dx.
\] 
Using this we have for $n \in \{1,\ldots,N-1\}$
\begin{eqnarray*}
\Pp \left[ \sum_i W_i =n \right] & \geq & \frac{N}{2} \frac{1}{N+1} {N \choose n} \int^1_0 
x^n(1-x)^{N-n} dx  \\
& = & \frac{N}{2} \frac{1}{N+1}{N \choose n} \frac{n! (N-n)!}{(N+1)!} \geq \frac{1}{4(N+1)}
\end{eqnarray*}
completing the proof of the claim. \qed

\noindent
\textbf{Remark 3.2.1} The proof uses only the special case of LNQD, that 
$\E[ X_i (\sum_{j\in A} X_j)^k] \leq 0$ for odd $k$ and $i \not \in J$.

\noindent 
\textbf{Remarks 3.2.2}
A converse inequality of the form $\E[[S]^p] \leq C_p \E[S^{2p}]$, which holds for i.i.d. variables, 
will certainly fail for negatively associated variables (as it does for the martingale case). Indeed 
the vector of i.i.d. variables, each taking values $\pm 1$ with probability $1/2$, 
conditioned that their sum equals $0$ will be negatively associated. 
This example is also useful to keep in mind if one wants to look for a maximal inequality
for negatively associated variables.
\subsection{Proof of Proposition \ref{NA}} \label{s3.3}
Following \cite{van2002randomly} we derive the negative association from the BKR inequality
first used in percolation models. 
and slight varaint on the types of negative dependency.
The BKR inequality \cite{reimer2000proof} is stated for a product probability measure $\mu$ on the product space
$\Omega = \prod_{i \in V} S_i $, where $V$, and $S_i$ for each $i \in V$, are finite sets. We recall the formulation:
for $\omega \in \Omega$ and $K \subseteq V$ the cylinder $[\omega]_K$ is defined by 
$[\omega]_K :=\{\omega': \omega'_i = \omega_i ,\, \forall i \in K\}$; for $A,B \subseteq \Omega$ the set
$A \Box B$ is defined by 
\[
A \Box B = \{\omega: \, \exists \;  \mbox{disjoint} \; K,L \subseteq V \, \mbox{with} \, 
[\omega]_K \subseteq A \; \mbox{and} [\omega]_L \subseteq B\}. 
\]
The BKR inequality states that
$ \mu(A \Box B) \leq \mu(A) \mu(B)$.

The papers \cite{kesten1998asymptotic,van2002randomly} define a product set-up as 
above within the graphical construction for their particle system. 
We argue slightly differently, 
by first defining a finite oriented percolation structure, equivalent to a system of
discrete time coalescing processes, for which BKR implies negative correlations. 
We then sketch the standard arguments for a sequence of discrete time approximations that converge
in distribution to the A particle process that we study. Since negative correlations 
statements survive such convergence we may deduce the inequalities we need. 

The product set-up we will construct will encode
a system of discrete time instantly coalescing particles performing 
random walks on $[-M,M]^d := \{-M, -M+1, \ldots, M-1,M\}^d$ over the time steps
$[0,T] := \{0,1,\ldots,T\}$. The particles will also have colours with values in
$[1,K] = \{1,2,\ldots,K\}$. 

In our construction, $V=V_1 \cup V_2$. The set $V_1 = [-M,M]^d \times [1,K]$ with $S_{(x,k)} = \{0,1\}$
for $(x,k)  \in V_1$ will encode the initial conditions; when $\omega_{(x,k)} = 1$ there
will be a particle of colour $k$ at position $x$ at time $0$. 
The set $V_2 = [0,T] \times [-M,M]^d \times 
[1,K]$ with $S_{(t,x,k)} = \{0,(\pm e_i)\} \times [1,K]$ for $(t,x,k) \in V_2$ will control the movement
and colour change of particles. Here $(e_1,\ldots,e_d)$ are the standard unit vectors.
 For example, when $\omega_{(t,x,k)} = (e_1,k')$ a particle
with colour $k$ at position $x$ at time $t$ will move to position $x+e_1$ at time $t+1$ and change colour to $k'$.

To define the positions of particles we adopt the terminology of oriented percolation. We first define 
paths running from time $0$ to time $T$:
a path $\pi$ consists of a sequence of space-colour points $\pi = \left((x_i,k_i): i = 0,1,\ldots,T\right)$ where $k_i \in [1,K]$ and
$x_i \in [-M,M]^d$ satisfy $|x_{i}-x_{i-1}| \in \{0,(\pm e_i)\}$ for $ i = 1,2,\ldots,T$.  The range of a path $R(\pi)$ is the set of time-space-colour points the path visits, that is $\{(i,x_i,k_i): i \in [0,T]\}$. For $\omega
\in \Omega$ we say a path $\pi = \left((x_i,k_i): i = 0,1,\ldots,T\right)$ is \textit{open} if
$\omega_{(x_0,k_0)} = 1 $ and 
\[
 (x_i - x_{i-1}, k_i)  = \omega_{(i-1,x_{i-1},k_{i-1})} \quad \mbox{for $i =1, \ldots,T$}.
\]
The key observation is that the open paths are instantly coalescing: if two open paths 
$\pi = \left((x_i,k_i): i = 0,1,\ldots,T\right)$ and $\pi' = \left((x'_i,k'_i): i = 0,1,\ldots,T\right)$ 
agree at time $t \in [0,T]$, that is $(x_t,k_t) = (x'_t,k'_t)$, 
then they agree at all times $t+1,\ldots, T$. 

Finally we define the occurence of a coloured particle
for $t \in [0,T],x \in [-M,M]^d,k \in [1,K]$ 
by
\[
\xi_t^{(k)}(x)(\omega) : = \I( \mbox{there is an open path $\pi$ with $(x_t,k_t)=(x,k)$}).
\]
The sum $\xi_t(x) : = \sum_{k=1}^K \xi_t^{(k)}(x)$ records the number of particles of any colour
at $x$ at time $t$.  
Note the variables $(\xi_t^{(k)}(x): x \in [-M,M], k \in [1,K])$ form a 
Markov chain over $t = 0,1,\ldots,T$. The initial conditions are
product Bernouilli variables since $\mu$ is a product measure. 
For suitably chosen $\mu$ we will show that
the process $(\xi_t(x): t \geq 0)$ will approximate our A-particle process. 

The following lemma is sufficient for all the negative dependency properties we use.
For a subset $A \subseteq [-M,M]^d \times [1,K]$ we define the set $O_A \subseteq \Omega$ by
\[
O_A = \{\omega: \mbox{for all $a \in A$ there are open paths $\pi_a$ whose final value 
$(x_T,k_T)$ equal $a$}\}.
\]

\begin{Lemma} \label{NDlemma}
Suppose $\mathcal{F}_1$ and $\mathcal{F}_2$ are two collections of subsets of 
$[-M,M] \times [1,K]$. Define
\[
\Omega_{1} = \bigcup_{A_1 \in \mathcal{F}_1}  O_{A_1}, \qquad \Omega_{2} = \bigcup_{A_2 \in \mathcal{F}_2}  O_{A_2}, \qquad \Omega_{1,2} = \bigcup_{\stackrel{A_1 \in \mathcal{F}_1, A_2 \in \mathcal{F}_2}{A_1 \cap A_2 = \emptyset}}  O_{A_1 \cup A_2}.
\]
Then $\mu(\Omega_{1,2}) \leq \mu(\Omega_{1}) \mu(\Omega_{2})$. 
\end{Lemma}

\textbf{Proof.} We claim $\Omega_{1,2} \subseteq \Omega_{1} \Box \Omega_{2}$, so that the result then follows from the BKR inequality. Indeed, 
for $\omega \in \Omega_{1,2}$ there exist disjoint $A_1,A_2$ so that
$\omega \in O_{A_1 \cup A_2}$. Let 
\[
R_i = \bigcup_{a \in A_i} R(\pi_a)
\]
be the union of the ranges of the paths $\pi_a$ that lead to $a \in A_i$, for $i = 1,2$. The coalescence of open paths implies that $R_1 \cap R_2 = \emptyset$. Moreover the 
cylinders $[\omega]_{R_i} \subseteq O_{A_i} \subseteq \Omega_{i}$ for $i =1,2$. \qed

As an example, one can apply this lemma to conclude, for $x \neq y$ and $m,n \geq 1$, that
\[
\Pp[ \xi_T(x) \geq m, \xi_T(y) \geq n] \leq \Pp[ \xi_T(x) \geq m] \, \Pp[\xi_T(y) \geq n] 
\]
by taking 
\begin{eqnarray*}
\mathcal{F}_1 &=& \{A = \{(k_1,x),\ldots, (k_m,x)\}: k_1, \ldots, k_m \in [1,K] \; \mbox{disjoint}\}, \\
\mathcal{F}_2 &=& \{A = \{(k_1,y),\ldots, (k_n,y)\}: k_1, \ldots, k_n \in [1,K] \; \mbox{disjoint}\},
\end{eqnarray*}
since then $\Omega_1 = \{ \xi_T(x) \geq m\}$ and $\Omega_2 = \{ \xi_T(y) \geq n\}$. 
To apply the lemma to conclude that
\[
\Pp[ \xi_T(x) \geq m_1 + m_2] \leq \Pp[ \xi_T(x) \geq m_1]  \, \Pp[\xi_T(x) \geq m_2] 
\]
we take $\mathcal{F}_i = \{A = \{(k_1,x),\ldots, (k_{m_i},x)\}: k_1, \ldots, k_{m_i} \in [1,K] \; \mbox{disjoint}\}$, noting then that
$\Omega_{1,2} = \{\xi_T(x) \geq m_1 + m_2\}$.

To apply the lemma to establish negative association for $\{\xi_T(x): x \in [-M,M]\}$ we 
fix disjoint subsets $\{i_1,\ldots,i_k\}$ and $\{j_1, \ldots,j_l\}$ of $\{1,\ldots,N\}$
and non-negative increasing functions $f:\R^k \to \R$ and $g:\R^l \to \R$.
We will show, for $m,n \geq 0$ 
 \begin{eqnarray*}
&& \hspace{-.4in}  \Pp[ f(\xi_T(x_{i_1}), \ldots, \xi_T(x_{i_k})) \geq m, \, g(\xi_T(x_{j_1}), \ldots, \xi_T(x_{j_l}))  \geq n] \\
& \leq & 
\Pp[ f(\xi_T(x_{i_1}), \ldots, \xi_T(x_{i_k})) \geq m] \, \Pp[ g(\xi_T(x_{j_1}), \ldots, \xi_T(x_{j_l}))  \geq n]
\end{eqnarray*}
which then implies the desired negative association. 
For a  set $A = \{(k_1,a_1),\ldots, (k_l,a_l)\} \subseteq [1,K] \times [-M,M]^d$ 
we associate a count 
\[
C^A(x) := \sum_{i=1}^{l} \I(x=a_i) \qquad \mbox{for $x \in [-M,M]^d$}.
\]
Then we choose
\begin{eqnarray*}
\mathcal{F}_1 &=& \{A = \{(k_1,a_1),\ldots, (k_l,a_l)\}: a_i \in \{x_{i_1},\ldots,x_{i_k}\} \; 
\mbox{for $i \leq l$,} \\
&& \hspace{.8in} (k_i,a_i)  \, \mbox{disjoint} \;  \mbox{and} \; f(C^A(x_{i_1}), \ldots, C^A(x_{i_k})) \geq m \}.
\end{eqnarray*} 
Since $f$ is increasing we find
$\Omega_1 = \bigcup_{A_1 \in \mathcal{F}_1} O_{A_1}$ is the event $\{f(\xi_T(x_{i_1}), \ldots, \xi_T(x_{i_k})) \geq m\}$. Use the equivalent definition for $\mathcal{F}_2$, replacing $f$, $m$ and  $\{x_{i_1},\ldots,x_{i_k}\}$ by 
$g$, $n$ and $\{x_{j_1},\ldots,x_{j_l}\}$. Note that $A_1 \in \mathcal{F}_1$
and $A_2 \in \mathcal{F}_2$ are automatically disjoint. 

We now sketch the steps needed to approximate our A-particle process by the 
discrete time processes above. As a summary, we go in the following order (i) pass to continuous time; (ii) 
check, by Dynkin's criterion, that $\sum_k \xi_t^{(k)}$ remains a Markov process; (iii) let
$K \to \infty$ to obtain the desired reaction and diffusion rates, and desired initial conditions; 
(iv) let $M \to \infty$ to obtain 
a system on $\Z^d$. 

The first aim (i) is to use a sequence of discrete time Markov chains, as in the oriented 
percolation construction above, to approximate a continuous time chain, 
which is easier to control when taking $K,M \to \infty$. 
We aim for the limit of the approximating sequence to be 
the continuous time Markov chain 
$(\bar{\xi}_t^{(k)}(x): x \in [-M,M]^d, k \in [1,K])_{0 \leq t \leq T}$,
with state space $\{0,1\}^{[-M,M]^d \times [1,K]}$; it will be a 
system of instantly coalescing coloured particles (where two particles at the 
same place and of the same colour instantly coalesce) 
performing random uniform colour updates at rate $\lambda$, and performing a simple random walk step 
with simultaneous random colour update at rate $2dD$. This can be defined by 
listing the rates of all types of possible jumps: independently for all $k,l \in [1,K]$ and for
$x, x+e_i \in [-M,M]^d$
\begin{eqnarray}
&& \!\!\left\{ 
\begin{array}{l} 
\bar{\xi}^{(k)}(x) \to \bar{\xi}^{(k)}(x) -1, \\
\bar{\xi}^{(l)}(x) \to (\bar{\xi}^{(l)}(x) +1) \wedge 1
\end{array} \right. \quad
\mbox{at rate $\lambda K^{-1} \bar{\xi}^{(k)}(x)$,} \nonumber \\
&&\!\! \left\{ 
\begin{array}{l} 
\bar{\xi}^{(k)}(x) \to \bar{\xi}^{(k)}(x) -1, \\
\bar{\xi}^{(l)}(x+e_i) \to (\bar{\xi}^{(l)}(x+e_i) +1) \wedge 1
\end{array} \right. \quad
\mbox{at rate $D K^{-1} \bar{\xi}^{(k)}(x)$.} \label{rates}
\end{eqnarray}

We use a standard method to approximate by discrete chains: choosing the time grid
$\{0,T/N, 2T/N, \ldots, T\}$ for the $N$'th approximation and linearly interpolation 
to produce a process indexed over $0 \leq t \leq T$. The $N$'th approximation has a 
transition matrix $P^{(N)}$ controlled by the choice of product measure $\mu^{(N)}$ in the percolation construction. 
We keep the initial conditions constant, that is product Bernoulli variables.  
Note the limit Markov chain has a 
finite state space and so it is enough to check that the elements $N P^{(N)}(\xi,\xi') \to Q(\xi,\xi')$
for the rates $Q(\xi,\xi')$ in (\ref{rates}) to ensure convergence of the processes (although we only use
convergence of the marginal at time $T$). 
The components of $\mu^{(N)}$ on the factor $S_{t,k,x}$ will have the form
$\mu_{(t,x,k)}((0,k)) = 1 - O(N^{-1})$, so that with high probability at most one of the 
open paths will either change location or change colour at any time. We may now choose
the values of $\mu_{(t,x,k)}((e_i,l))$ and $\mu_{(t,x,k)}((k,l))$, both $O(N^{-1})$, to yield the desired 
limiting $Q$ matrix.  

In step (ii) we apply Dynkin's criterion \cite{dynkin1965markov} to see that the process 
$\bar{\xi}(x) = \sum_{k=1}^K \bar{\xi}^{(k)}(x) $
is a Markov process with state space $\{0,1,\ldots,K\}^{[-M,M]^d}$ and with jump rates
\begin{eqnarray*}
&& \!\!\left\{ 
\begin{array}{l} 
\bar{\xi}(x) \to \bar{\xi}(x) -1, \\
\bar{\xi}(x+e_i) \to \bar{\xi}(x+e_i) +1,
\end{array} \right. \quad
\mbox{at rate $D K^{-1} \bar{\xi}(x)(K-\bar{\xi}(x+e_i))$,} \\
&& 
\bar{\xi}(x) \to \bar{\xi}(x) -1
 \quad
\mbox{at rate $\lambda K^{-1} \bar{\xi}(x)(\bar{\xi}(x)-1)) +
D K^{-1} \bar{\xi}(x) \sum_{y \sim x: y \in [-M,M]^d} \bar{\xi}(y)$.}
\end{eqnarray*}
The initial conditions  for $\{\bar{\xi}_0(x): x\in [-M,M]^d\}$ are now independent Poisson-Bernouilli variables. 

Step (iii) is to choose $D = \frac{D_A}{2d}$ and $\lambda = K \lambda_A$ and let $K \to \infty$.
This yields a Markov process with (infinite) state space $\{0,1,\ldots,\}^{[-M,M]^d}$ and jump rates
\begin{eqnarray*}
&& \!\!\left\{ 
\begin{array}{l} 
\bar{\xi}(x) \to \bar{\xi}(x) -1, \\
\bar{\xi}(x+e_i) \to \bar{\xi}(x+e_i) +1,
\end{array} \right. \quad
\mbox{at rate $\frac{D_A}{2d} \, \bar{\xi}(x)$,} \\
&& 
\bar{\xi}(x) \to \bar{\xi}(x) -1
 \quad
\mbox{at rate $\lambda_A \bar{\xi}(x)(\bar{\xi}(x)-1))$.}
\end{eqnarray*}
The limiting initial conditions $\xi_0(x)$ will be independent for $x \in [-M,M]^d$. We choose them to be 
identically distributed as a law lying in the closure $\overline{\mathcal{M}_{PB}}$ of the set of Poisson-Bernouilli laws. 
Note that the the total number of particles at $t=0$ is almost surely finite 
and the total number then only decreases in time due to coalescence.  
This is the desired A-particle process restricted to the region $[-M,M]^d$ (that is jumps
leaving this region are suppressed). 

The final step (iv), letting $M \to \infty$ to yield the process on $\Z^d$ can be done, for example, by coupling
strong equations, as in (\ref{xieqn}), for a finite $(\bar{\xi}_M(t,x):x \in [-M,M]^d, t \geq 0)$ to
an infinite system $(\bar{\xi}(t,x):x \in \Z^d, t \geq 0)$, using the same 
Poisson clocks where possible,  and estimating the difference $\E[ \|\bar{\xi}_M(t) - \bar{\xi}(t)\|_{\alpha} ] $ at later times in
this weighted norm. Such a technique is well known when approximating infinite 
volume spatial systems and we omit details (which are in the thesis \cite{lukins2017thesis}
for the case $d=2$.)

For each of the steps above that establish a convergence of processes, the negative correlation
statements can be carried over, leading to the conclusions of Proposition \ref{NA}.
\section{Proof of Theorem \ref{mainresult}}  \label{s4}

We suppose throughout this section that $(\xi,\eta)$ are solutions to (\ref{xieqn}) and (\ref{etaeqn}) in $d \geq 3$ satisfying the hypotheses of Theorem \ref{mainresult}.  
We need to develop over time various expectations of moments. These are infinitesimal generator calculations 
which we do by calculus from the strong equations. 
Compensating each Poisson process in (\ref{xieqn}) and (\ref{etaeqn}), that is writing $P_t$ as $(P_t- \lambda t) +\lambda t$ 
where $\lambda$ is the corresponding intensity, we can define (local) martingales $dM^A_t(x),dM^B_t(x)$ so that 
\begin{eqnarray}
d \xi_t(x) & = & D_A \Delta \xi_t(x)dt - \lambda_A \xi_t(x)(\xi_t(x)-1) dt + dM^A_t(x)  \label{mgpa}\\
d \eta_t(x) & = & D_B \Delta \eta_t(x)dt - \lambda_B \xi_t(x) \eta_t(x) dt + dM^B_t(x) \label{mgpb}
\end{eqnarray}
where $\Delta:=\Delta_x$ is the discrete Laplacian defined by $ \Delta f(x) = (1/2d) \sum_{y: y \sim x} (f(y) - f(x))$.
We will need the brackets processes defined by 
\[
[ \xi_t(x),\xi_t(y)] = \sum_{s \leq t} (\xi_t(x)- \xi_{t-}(x))(\xi_t(y) - \xi_{t-}(y)).
\]
The (unique continuous adapted) compensators  $\llangle \xi(x), \xi(y) \rrangle_t$ 
can be calculated to be 
\begin{equation} 
d \llangle \xi(x), \xi(y) \rrangle_t = \left\{ \begin{array}{ll}
\frac{D_A}{2d} \sum_{y: y \sim x} (\xi_t(y)+ \xi_t(x)) dt + \lambda_A \xi_t(x) (\xi_t(x)-1) \, dt & \mbox{if $x=y$,} \\
 - \frac{D_A}{2d} \left( \xi_t(x) + \xi_t(y) \right) \, dt & \mbox{if $x \sim y$,} \\
0 & \mbox{otherwise.}
\end{array} \right. \label{xicomp}
\end{equation}
We will not need $[\eta_t(x),\eta_t(y)]$ and since
there are no simultaneous jumps between $\xi_t(x)$ and $\eta_t(y)$ we have $[\xi_t(x),\eta_t(y)] = 0$. 

\noindent
\textbf{Notation.} To ease the presentation we use the following notation: 

We write $\hat{\xi}_t$ and $\hat{\eta}_t$ for the expected values $\E[\xi_t(x)]$ and $\E[\eta_t(x)]$.  

For bounded $f,g:\Z^d \to \R$ we write $\La f,g \Ra$ for the sum $\sum_x f(x) g(x)$ when this is defined. 

When $\phi: \Z^d \times \Z^d \to \R$ we similarly write $\La f * g, \phi \Ra$ for the sum
$\sum_{x,y} f(x) g(x) \phi(x,y)$.
\subsection{$A$ particle estimates} \label{s4.1}
As in Bramson and Griffeath \cite{bramson1980asymptotics}, we start with a 'crude' upper 
and lower bound on $\hat{\xi}_t$: there exist 
$0 < c_1= c_1(\lambda_A, D_A, \hat{\xi}_0)$ and $0 < c_2= c_2(\lambda_A,D_A,\mathcal{L}(\xi_0))$ so that
\begin{equation} \label{crudeLUB}
\frac{c_1}{t} \leq \hat{\xi_t} \leq \frac{c_2}{t} \quad \mbox{for all $t \geq 1$.}
\end{equation}
These bounds already shows that the correct polynomial rate is $O(t^{-1})$ and constituted a 
key part of the argument in  \cite{bramson1980asymptotics}. The
efforts in the modified rate equation method are to resolve the exact constant in the asymptotics. 

Although the bounds in \cite{bramson1980asymptotics} are derived for the instantaneously coalescing 
case, they were extended in \cite{kesten1998asymptotic} to hold for the models considered there and can also
be adapted for our $A$ particle model. Indeed the lower bound
follows immediately from negative correlation, as we will see just below. The upper bound argument
in \cite{bramson1980asymptotics} is however rather closely related to the methods here; they show that
the expected density satisfies a discrete time version of the equation 
$d \hat{\xi_t} \leq - C \hat{\xi_t}^2$. Regardless of the value of the constant $C>0$, this already implies an upper bound as in
(\ref{crudeLUB}). It is not surprising therefore that the related arguments in this paper can also be used to derive such an upper bound. We show the steps needed in the Appendix, mainly so 
that the paper can be considered self-contained (rather than suggesting that the steps here are 
any easier than adjusting the arguments from \cite{bramson1980asymptotics}). 

We suppose now that the bounds
(\ref{crudeLUB}) have been established. 
As a consequence, using the negative dependency estimate (\ref{natwo}) we find constants $c_3 = c_3(n,\lambda_A,D_A,\mathcal{L}(\xi_0))$ so that, 
for $x \in \Z^d, \, n \in N$ and $t \geq 1$,
\begin{equation} \label{moments}
\E \left[ \xi_t(x) (\xi_t(x)-1) \ldots (\xi_t(x)-(n-1))\right] \leq \frac{c_3}{t^n} \quad \mbox{and} \quad
\E \left[ \xi_t^n(x)\right] \leq \frac{c_3}{t}.
\end{equation}
Here the second bound follows by expressing the moments in terms of factorial moments.
%

\noindent
\textbf{One point estimates.}
Using the moments assumptions on our initial conditions, the local martingales in (\ref{mgpa}) and (\ref{mgpb}) 
are true martingales. In particular, taking expectations in (\ref{mgpa}) and exploiting translation invariance, we have
\begin{equation} \label{basicidea1}
d \hat{\xi}_t = - \lambda_A \E[ \xi_t(0) (\xi_t(0)-1)] dt.
\end{equation} 
In particular $t \to \hat{\xi}_t$ is decreasing, and using negative dependence (\ref{natwo}) we obtain
\[
\frac{d \hat{\xi}_t}{dt}  \geq - \lambda_A \hat{\xi}_t^2 \quad \mbox{for $t \geq 0$}
\]
and solving this differential inequality yields 
\begin{equation} \label{LB}
\hat{\xi_t} \geq \left( \hat{\xi}_0^{-1} + \lambda_A t \right)^{-1}
\end{equation}
(which is one way to establish the lower bound in (\ref{crudeLUB})).

For a (suitably smooth and integrable) test function $\phi_t(x)$ we have 
\begin{equation} \label{onepointA}
d \La \xi_t,\phi_t \Ra  = \La \xi_t, \dot{\phi}_t + D_A \Delta \phi_t \Ra dt - \lambda_A
\La \xi_t (\xi_t-1), \phi_t \Ra dt + d \La M^A_t,\phi_t \Ra.
\end{equation}
A useful example is to take integrable $f:\Z^d \to \R$ and choose
\[
\phi_{s}(x) = P^A_{t-s}f(x) := \sum_y p^A_{t-s}(y) f(x-y) \quad \mbox{for $s \in [0,t]$}
\]
where $p^A_{t}(y)$ is the transition density for an $A$ particle performing rate $D_A$ simple random walk. 
Then $\dot{\phi}_s + D_A \Delta \phi_s=0$ over $s \in [0,t]$ so that
\begin{equation} \label{GreenA}
d \La \xi_s,P^A_{t-s}f \Ra = - \lambda_A \La \xi_s (\xi_s-1), P^A_{t-s}f \Ra ds
 + d \La M^A_s,P^A_{t-s}f \Ra.
\end{equation}
Note the drift term on the right hand side is non-positive when $f \geq 0$. 
This leads to the following one point estimate, which allows us, 
over short periods of time, to approximate the 
solution by (random walk) heat flow, bounding the error caused by coalescence. 
\begin{Lemma} \label{A1pt}
For $f \geq 0$ we have $\E \left[ \La \xi_t,f \Ra \right] \leq  \E \left[ \La \xi_{t-s}, P^A_s f \Ra \right] $. Moreover
\begin{equation} \label{29}
\left| \E \left[ \La \xi_t,f \Ra \right] - \E \left[ \La \xi_{t-s}, P^A_s f \Ra \right] \right| 
\leq 4  c_2^2 \lambda_A \langle f,1 \rangle s t^{-2} \quad \mbox{when for $0<s < t/2$, $t \geq 1$.} 
\end{equation}
\end{Lemma}
\noindent
\textbf{Proof} 
 Integrating (\ref{GreenA}) over the interval $[t-s,t]$ and taking expectations we reach
\[
 \E \left[ \La \xi_t,f \Ra \right] - \E \left[ \La \xi_{t-s}, P^A_s f \Ra \right] 
= - \lambda_A
\int^t_{t-s} \E \left[ \La \xi_r (\xi_r-1), P^A_{t-r}f \Ra \right] dr \leq 0.
\]
Using translation invariance and the negative dependence result (\ref{natwo}) we find
\[
 \E \left[ \La \xi_r (\xi_r-1), P^A_{t-r}f \Ra \right] =  \E \left[ \xi_r(0) (\xi_r(0)-1) \right] \La 1,  P^A_{t-r}f \Ra
 \leq \hat{\xi}^2_r \, \La 1, f \Ra.
\]
The inequality in (\ref{29}) follows from the upper bound (\ref{crudeLUB}). 
 \qed
%

\noindent
\textbf{Two point estimates.} 
We aim to use a two-point calculation to estimate 
$\E[ \xi_t(x)(\xi_t(x)-1)]$. Use the integration by parts
$d(\xi_t(x) \xi_t(y)) = \xi_{t-}(x) d\xi_t(y) + \xi_{t-}(y) d\xi_t(x) + d[\xi(x),\xi(y)]_t$
to find
\begin{eqnarray*}
d( \xi_t(x) \xi_t(y)) & = & D_A \Delta (\xi_t(x) \xi_t(y)) dt - \lambda_A \xi_t(x) \xi_t(y)
(\xi_t(x) + \xi_t(y) -2) dt \\
& &  \;\; + \xi_{t-}(x) dM^A_t(y) + \xi_{t-}(y) dM^A_t(x)  + d[\xi(x),\xi(y)]_t
\end{eqnarray*}
where the Laplacian $\Delta=\Delta_x + \Delta_y$ now acts in both variables. 
We replace the jump term 
$d[\xi(x),\xi(y)]_t$ by its compensator $d\llangle\xi(x),\xi(y)\rrangle_t$ as in (\ref{xicomp}) plus martingale increments, which 
we henceforth denote by $MI$, reaching
\begin{equation} \label{dxixy}
d( \xi_t(x) \xi_t(y))  =  D_A \Delta (\xi_t(x) \xi_t(y)) dt - \lambda_A \xi_t(x) \xi_t(y)
(\xi_t(x) + \xi_t(y) -2) dt + d \llangle \xi(x), \xi(y) \rrangle_t + MI.
\end{equation}
Using also discrete integration by parts to transfer the discrete Laplacian onto 
the test function, that for (suitably smooth and integrable) $\phi_t(x,y)$ 
\begin{eqnarray}
d \La \xi_t * \xi_t,\phi_t \Ra & = &  \La \xi_t * \xi_t, \dot{\phi}_t + D_A \Delta \phi_t \Ra dt 
- \lambda_A \sum_x \sum_y \xi_t(x) \xi_t(y) (\xi_t(x)+\xi_t(y)-2) \phi_t(x,y) dt \nonumber \\
& & \;\; + \sum_x \left( \frac{D_A}{2d} \sum_{y:y \sim x} (\xi_t(y) + \xi_t(x))
 + \lambda_A \xi_t(\xi_t(x)-1)\right)  \phi_t(x,x) dt
\nonumber \\
&&  \;\;\;\; - \frac{D_A}{2d} \sum_x \sum_{y:y \sim x} (\xi_t(x) + \xi_t(y)) \phi_t(x,y)dt + \mbox{MI}
\nonumber \\
& = & \La \xi_t * \xi_t, \dot{\phi}_t + D_A \Delta \phi_t \Ra dt 
 - \lambda_A \sum_x \sum_{y: y \neq x} \xi_t(x) \xi_t(y) (\xi_t(x)+\xi_t(y)-2) \phi_t(x,y) dt \nonumber \\
& & \;\; - \lambda_A \sum_x \left( 2\xi^3_t(x) -  3\xi^2_t(x) +  \xi_t(x) \right) \phi_t(x,x) 
 + D_A \La \xi_t,\Pi \phi_t \Ra dt + \mbox{MI} \label{39a}
\end{eqnarray}
where we define
\[
\Pi \phi_t(x) := \frac{1}{2d} \sum_{y:y \sim x} \left( \phi_t(x,x) + \phi_t(y,y) - \phi_t(x,y) - \phi_t(y,x) \right).
\]
A useful specific choice of test function is as follows: set 
\begin{equation} \label{psidefn}
\psi^{a,b}_t(x,y) = \Pp[S^{x,A}_t =a, \, S^{y,A}_t = b, \, \tau >t]
\end{equation}
where $S^{x,A}$ and $S^{y,A}$ are independent type $A$ random walks started at $x$ and $y$, and where
\[
\tau = \inf\{t: \int^t_0 I(S^x_s=S^y_s) ds > \mathcal{E}\}
\]
for $\mathcal{E}$ is an independent exponential variable with rate $2 \lambda_A$.
Then choosing $\phi_s=\psi^{a,b}_{t-s}$ over $s \in [0,t]$ we have that
\begin{equation} \label{phieqn}
\dot{\phi}_s(x,y)  + D_A \Delta \phi_s(x,y) = 2 \lambda_A \phi_s(x,y) I(x=y), \quad \mbox{and $\phi_t(x,y) = I(x=a, y=b)$.}
\end{equation}
Using this test function in (\ref{39a}) we find, for $s \in [0,t]$,
\begin{eqnarray}
d \La \xi_s * \xi_s,\phi_s \Ra & = & 
- \lambda_A \sum_x \sum_{y:y \neq x} \xi_s(x) \xi_s(y) (\xi_s(x)+\xi_s(y)-2) \phi_s(x,y) ds 
\nonumber \\
&& \hspace{-.3in}- \lambda_A \sum_x \left(2\xi^3_s(x) - 5 \xi_s^2(x) + \xi_s(x) \right) \phi_s(x,x) ds
 + D_A \La \xi_s, \Pi \phi_s \Ra ds + \mbox{MI.} \label{39b}
\end{eqnarray}
We now define $\tilde{\phi}:[0,t] \times \Z^d \to \R$ by $\tilde{\phi}_s(x)= \phi_s(x,x)$. The idea is to choose $a=b=0$
so that 
\begin{equation} \label{endpoint}
\E[ \La \xi_t * \xi_t,\phi_t \Ra -  \La \xi_t,\tilde{\phi}_t \Ra ] = E[ \xi_t(0)(\xi_t(0)-1)].
\end{equation}
A short calculation shows that
\[
\Delta \tilde{\phi}_s(x) = \Delta \phi_s(x,x) + \Pi \phi_s(x).
\]
Then $\dot{\tilde{\phi}}  + D_A \Delta \tilde{\phi} = 2 \lambda_A \tilde{\phi} + D_A \Pi \phi$ and
using (\ref{onepointA}) we have
\begin{equation} \label{39c}
d \La \xi_s,\tilde{\phi}_s \Ra  = 2 \lambda_A \La \xi_s,\tilde{\phi}_s \Ra ds - \lambda_A
\La \xi_s(\xi_s-1),\tilde{\phi}_s \Ra  ds + D_A \La \xi_s,\Pi \phi_s \Ra ds + \mbox{MI}.
\end{equation}
Combining (\ref{39b}) and (\ref{39c}) we find
\begin{eqnarray} \label{A2pteqn}
 d \La \xi_s * \xi_s,\phi_s \Ra - d \La \xi_s,\tilde{\phi}_s \Ra & = & 
- \lambda_A \sum_x \sum_{y:y \neq x} \xi_s(x) \xi_s(y) (\xi_s(x)+\xi_s(y)-2) \phi_s(x,y) ds \nonumber \\ 
&&  \;\; - 2 \lambda_A  \La \xi_s(\xi_s-1)(\xi_s-2),\tilde{\phi}_s \Ra ds + \mbox{MI} \label{39end}
\end{eqnarray}
Note that both drift terms on the right hand side are non-positive.
This will lead to the following two point moment estimate.
%
\begin{Lemma} \label{A2pt}
Let $\psi^{0,0} $ be defined as in (\ref{psidefn}). Then 
\begin{equation} \label{A2ptUB}
\E \left [\xi_t(0) (\xi_t(0)-1) \right] \leq \E \left[ \La \xi_{t-s} * \xi_{t-s},\psi^{0,0}_s \Ra \right]
\end{equation}
and moreover, for some $c_4 = c_4(\lambda_A,D_A, \mathcal{L}(\xi_0))$, for $0<s<t/2$ and $t \geq 1$ 
\begin{equation} \label{A2ptestimate}
\left| \E \left [\xi_t(0) (\xi_t(0)-1) \right] - \E \left[ \La \xi_{t-s} * \xi_{t-s},\psi^{0,0}_s \Ra \right] \right| \leq 
c_4 \left( t^{-1} s^{-d/2}+  s \, t^{-3} \right). 
\end{equation}
\end{Lemma}
\noindent
\textbf{Proof.} Choosing $\phi_s=\psi^{0,0}_{t-s}$ as above,
integrating (\ref{A2pteqn}) over the interval $[t-s,t]$ and then take expectations we reach, using (\ref{endpoint}),
\begin{eqnarray}
&& \E \left [\xi_t(0) (\xi_t(0)-1) \right]  - \E \left[ \La \xi_{t-s} * \xi_{t-s},\psi^{0,0}_s \Ra \right] \nonumber \\
& = &
- \E \left [\La \xi_{t-s}, \tilde{\phi}_{t-s} \Ra \right] 
 - \lambda_A \int^t_{t-s}  \sum_x \sum_{y: y \neq x}  \E \left [ \xi_r(x) \xi_r(y) (\xi_r(x)+\xi_r(y)-2) \right]  \phi_r(x,y)  dr \nonumber \\
&& \hspace{.3in} - 2  \lambda_A \int^t_{t-s} \E \left [\La \xi_r(\xi_r-1)(\xi_r-2),\tilde{\phi}_r \Ra \right] \leq 0. \label{temp107}
\end{eqnarray}
Now we use $\psi^{0,0}_t(x,y) \leq p^A_t(x) p^A_t(y)$ and the bound
\begin{equation} \label{pAbound}
p^A_t(x) \leq C(D_A) (t^{-d/2} \wedge 1) \qquad \mbox{for all $t,x,y$.}
\end{equation} 
In particular $\La 1, \tilde{\phi}_{t-s} \Ra = \sum_x \psi^{0,0}_{s}(x,x) \leq C(D_A) s^{-d/2}$. The negative association and 
(\ref{natwo}) imply, when  $x \neq y$, that
\begin{equation} \label{41.5}
 \E \left [ \xi_r(x) \xi_r(y) (\xi_r(x)+\xi_r(y)-2) \right]  =  \E \left [ \xi_r(x)(\xi_r(x)-1) \xi_r(y)  \right]
 + \E \left [ \xi_r(x) \xi_r(y) (\xi_r(y)-1) \right] \leq 4 \hat{\xi}^3_r.
 \end{equation}
Now using (\ref{moments}) all three terms on the right hand side of (\ref{temp107}) can thus be bounded as stated.  
\qed

The final estimate is on the non-coalescence probability test function used above, namely
$\psi^{0,0}_t(x,y) = \Pp[S^{x,A}_t =0, \, S^{y,A}_t = 0, \, \tau >t]$. 
The random walk arguments from \cite{kesten1998asymptotic} give good 
approximations for this $\psi$; it is shown, when $d \geq 3$, that for some $\delta>0$ and $c_5 = c_5(D_A)$
\begin{equation} \label{KdBl12}
\sum_{x,y \in\Z^d} \left| \psi^{0,0}_t(x,y) - p_A p^A_t(x) p^A_t(y) \right| \leq c_5 \, t^{-\delta} \quad 
\mbox{for all $t \geq 0$}.
\end{equation}

\noindent
\textbf{Remark. 4.1.1} Indeed, taking $m=0$ in Lemma 12 in \cite{kesten1998asymptotic} shows that
$\delta = (d-2)/(3d^2-3d-4)$ will work. However \cite{kesten1998asymptotic} treats a more general random walk jump 
requiring only second moments.  In our nearest neighbour case the argument yields $c_5=c_5(D_A,\delta)$ for any 
$0<\delta < \frac12 \frac{(d-2)}{(d-1)}$.

\subsection{$A$-particle modified rate equation} \label{s4.2}
We repeat the lines of the argument of van den Berg and Kesten \cite{kesten1998asymptotic}, 
which will help prepare us for related steps when deriving the $B$-particle modified rate equation.
There are four approximations:
\begin{eqnarray}
\frac{d \hat{\xi}_t}{dt} & = &  - \lambda_A \E[\xi_t(0)(\xi_t(0)-1)] \nonumber \\
& \approx & - \lambda_A \E[\La \xi_{t-s} * \xi_{t-s}, \psi^{0,0}_s \Ra]  \label{app1} \\
& \approx & - \lambda_A p_A  \E[ \La \xi_{t-s},p^A_s \Ra^2]  \label{app2} \\
& \approx & - \lambda_A p_A  (\E[ \La \xi_{t-s},p^A_s \Ra])^2   \label{app3} \\
&   = &- \lambda_A p_A  \hat{\xi}_{t-s}^2 \nonumber \\
&  \approx & - \lambda_A p_A  \hat{\xi}_{t}^2 \label{app4} 
\end{eqnarray}
Using the lemmata from the previous section we will below bound the errors in these approximations to show, 
for a $\kappa_1 \in (0,1)$ to be chosen soon, that
\begin{equation}
\label{Aere}
\frac{d \hat{\xi}_t}{dt} =  - \lambda_A p_A  \hat{\xi}_{t}^2 + \mathcal{E}_t, \quad \mbox{where $\mathcal{E}_t = O(t^{-2 - \kappa_1})$.}
\end{equation}
Now a short calculus exercise (as van den Berg and Kesten put it), using 
the lower bound $\hat{\xi}_t \geq c_1/t$, shows that
\[
\left| \hat{\xi}_t - \frac{1}{\lambda_A p_A t} \right| \leq C t^{-1-\kappa_1}
\]
completing the proof of the asymptotics for the $A$-particle density. One way to do this calculus exercise 
is to 
 integrate $\hat{\xi}_t^{-2} d \hat{\xi}_t/dt$ to get 
\begin{equation} \label{neededA}
\hat{\xi}_{t}^{-1} - \hat{\xi}_{t_0}^{-1} = \lambda_A p_A (t-t_0) - \int^t_{t_0} \hat{\xi}_s^{-2} \mathcal{E}_s ds
\end{equation}
and use the lower bound to see that the final integral is bounded by $C t^{1-\kappa_1}$.

While checking (\ref{Aere}) we shall take $t \geq 1$ and make use of a chosen $0<s \leq t/2$. 
We will use a running constant $C$ which will
depend on $d, \lambda_A,D_A, \mathcal{L}(\xi_0)$. 
Lemma \ref{A2pt} gives the error bound for the first approximation (\ref{app1}), namely
$C \left( t^{-1} s^{-d/2}+  s \, t^{-3} \right)$. 
Using (\ref{KdBl12}), the negative association, and the bound $\psi^{0,0}_t(x,y) \leq (p^A_s(x))^2$  we bound the approximation in (\ref{app2}) by
\begin{eqnarray*}
\left| \, \E[\La \xi_{t-s} * \xi_{t-s}, \psi^{0,0}_s \Ra]  -   p_A \E[ \La \xi_{t-s},p^A_s \Ra^2] \, \right|
& \leq &  c_4 s^{-\delta} \hat{\xi}_{t-s}^2 + \sum_x \E[\xi_{t-s}^2(x)] (p^A_s(x))^2 \\
& \leq & C (t^{-2} s^{-\delta} + t^{-1} s^{-d/2})
\end{eqnarray*}
where we have used (\ref{moments}) and (\ref{pAbound}) for the final inequality.
%

The negative association leads to a simple 
estimate of the variance of $ \La \xi_t,f \Ra$. Indeed we have for $t \geq 1$ and $f \geq 0$ 
\begin{eqnarray}
\mbox{Variance}( \La \xi_t,f \Ra) 
&=&  \sum_x \sum_y \E[(\xi_t(x) - \hat{\xi}_t) (\xi_t(y) - \hat{\xi}_t)] f(x) f(y) \nonumber \\
& \leq & \sum_x \E[(\xi_t(x) - \hat{\xi}_t)^2] f^2(x) \nonumber \\
& \leq & \E[(\xi_t(0))^2]   \La f^2,1 \Ra.
 \label{Avariance}
\end{eqnarray}
Now we use (\ref{Avariance}) to bound the third approximation in (\ref{app3}) by
\begin{equation}
\left| \E \left[ \La \xi_{t-s},p^A_s \Ra^2 \right] - 
\E \left[ \La \xi_{t-s},p^A_s \Ra \right]^2 \right|
\leq   \hat{\xi}_{t-s} \La 1, (p^A_s)^2 \Ra \leq C t^{-1} s^{-d/2}. \label{Astep3}
\end{equation} 
Using translation invariance, Lemma \ref{A1pt} implies that
$ | \hat{\xi}_t -  \hat{\xi}_{t-s}| \leq   4  c_2^2 \lambda_A s t^{-2}$.
Hence we may bound the final approximation in (\ref{app4}) using
\[
\left| \hat{\xi}_t^2 -  \hat{\xi}_{t-s}^2 \right|
= \left| \hat{\xi}_t -  \hat{\xi}_{t-s} \right| \left| \hat{\xi}_t +  \hat{\xi}_{t-s} \right| \leq C s t^{-3}.
\] 
Combining all four estimates we reach (\ref{Aere}) with
 $| \mathcal{E}_t| \leq C ( t^{-1} s^{-d/2} + s t^{-3} + t^{-2}s^{-\delta})$. 
Choosing $s=t^{\alpha}$ where $\alpha \in (2/d,1)$ we have
\begin{equation} \label{finalEA}
\left| \mathcal{E}_t \right| \leq C t^{-2-\kappa_1} \quad 
\mbox{where $\kappa_1 = \min \left\{\frac{\alpha d}{2}-1, 1-\alpha,  \delta \alpha\right\}>0$.} 
\end{equation}

\noindent
\textbf{Remark 4.2.1} Choosing $\delta$ just below $\frac12 \frac{d-2}{d-1}$, as in the remark after (\ref{KdBl12}), we can then choose an optimal value of $\alpha$ above (just above $\frac{2d-2}{3d-4}$) to find $\kappa_1$ as close to $\frac{d-2}{3d-4}$ as we want.

%
%
%
\subsection{$B$-particle estimates} \label{s4.3} 

A key estimate in this paper is to bound the covariance between $\La \xi_t,f \Ra$ and $ \La \eta_t,g \Ra$.
Without negative dependence results for $B$ particles we rely simply on H\"{o}lder's inequality. 

Recall we have translation invariant solutions and we are writing $\hat{\eta}_t$ and $\hat{\xi}_t$ for $E[\eta_t(0)]$
and  $E[\xi_t(0)]$.
Suppose $f,g: \Z^d \to [0,\infty)$ satisfy $\La f,1\Ra =  \La g,1\Ra=1$. Choose $p,q$ satisfying $p^{-1} + q^{-1} =1$.
Then   
\begin{eqnarray}
\left| \E \left[ \La \xi_t,f \Ra \La \eta_t,g \Ra \right] - \E \left[ \La \xi_t,f \Ra \right] \, 
\E\left[ \La \eta_t,g \Ra \right] \right| 
& = & \left| \E\left[ \La \xi_t- \hat{\xi}_t,f \Ra  \La \eta_t,g \Ra \right] \right|  \nonumber \\
& \leq & \left| \E\left[ \La \xi_t - \hat{\xi}_t,f \Ra ^q \right]  \right|^{1/q} \, \left| 
\E\left[ \La \eta_t,g \Ra^p \right] \right|^{1/p}. \label{StepA}
\end{eqnarray}
Let $\tilde{\eta}= \E[ \eta | \sigma(\xi)]$ be the $B$ particle process conditional on the 
$A$ population. 
In the above, if we take conditional expectations with respect to $\sigma(\xi)$ before applying
H\"{o}lder's inequality then we may replace $\La \eta,g \Ra$ by $ \La \tilde{\eta},g \Ra$
on the right hand side of (\ref{StepA}). 
Taking conditional expectations in (\ref{mgpb}) we see that
$\tilde{\eta}$ solves 
\[
d \tilde{\eta}_t(x) = D_B \Delta \tilde{\eta}_t(x)dt - \lambda_B \xi_t(x) \tilde{\eta}_t(x) dt 
\]
namely the heat flow with random killing at rate $\lambda_B \xi_t(x)$ and constant initial condition $\E[\eta_0(0)]$.
In particular $\tilde{\eta}_t(x) \leq \E[\eta_0(0)]$ for all $t,x$. Then
\begin{eqnarray}
\E \left[ \La \tilde{\eta}_t,g \Ra ^p \right] 
& \leq & \E \left[  \La \tilde{\eta}_t^p,g \Ra  \right]  \nonumber  \\ 
& \leq & \E \left[  \La \tilde{\eta}_t,g \Ra \right] \E[\eta_0(0)]^{p-1}  \nonumber  \\
& = &  \E \left[  \La \eta_t,g \Ra \right] \E[\eta_0(0)]^{p-1} = \E[\eta_0(0)]^{p-1} \, \hat{\eta}_t. \label{StepB}
\end{eqnarray}
The idea is that the $L^p$ norm $\| \La \tilde{\eta},g  \Ra \|_p$ will be bounded by $C \hat{\eta}_t^{1/p}$. By taking 
$p$ close to $1$ we will get as close as desired to polynomial rate of $t^{-\theta}$, which will be needed
if the error estimates that use H\"{o}lder's inequality are small for the modified $B$ particle rate equation.
For this reason we will need high $q$'th moments of the $A$ population, which is why we asked for all
moments of the initial conditions to be finite. 

To bound the $L^q$ norm of $ \La \xi_t- \hat{\xi}_t,f \Ra$ 
we can use the square function inequality Proposition \ref{MZ} from Section \ref{s3.2} 
to bound, when $q/2 \in \N$,
\begin{eqnarray}
\E \left[  \La \xi_t - \hat{\xi}_t,f \Ra^q \right] & \leq &
 C_q \E \left[ \La(\xi_t - \hat{\xi}_t)^2, f^2 \Ra^{q/2} \right] \nonumber \\
 & \leq & C_q \E \left[ \La \xi_t^2, f^2 \Ra^{q/2} \right] + C_q t^{-q} \La 1, f^2 \Ra^{q/2}   \nonumber  \\
 & \leq & C_q \E [ \xi_t^q(0)]  \La 1,f^2 \Ra^{q/2} + C_q t^{-q} \La 1, f^2 \Ra^{q/2} \nonumber \\
 & \leq & C_q t^{-1} \La 1,f^2 \Ra^{q/2} + C_q t^{-q} \La 1, f^2 \Ra^{q/2}   \label{oneit}
\end{eqnarray}
by the moment bounds (\ref{moments}). Or we can do better by repeating the trick: 
write $\xi_t^2(x) = \xi_t^2(x) - m_2(t) + m_2(t) $ where $m_2 (t) := \E[\xi^2_t(0)] = O(t^{-1})$
and use the square function inequality on the negatively associated variables $(\xi_t^2(x) - m_2(t): x \in \Z^d)$. 
This gives, when $q/4 \in \N$, 
\begin{eqnarray}
\E \left[  \La \xi_t - \hat{\xi}_t,f \Ra^q \right] & \leq &
 C_q \E \left[ \La \xi_t^2 - m_2(t), f^2 \Ra^{q/2} \right] + C_q t^{-q/2} \La 1, f^2 \Ra^{q/2}  \nonumber  \\
 & \leq &
 C_q \E \left[ \La(\xi_t^2 - m_2(t))^2, f^4 \Ra^{q/4} \right] + C_q t^{-q/2} \La 1, f^2 \Ra^{q/2}   \nonumber \\
 & \leq &
 C_q \E \left[ \La(\xi_t^4, f^4 \Ra^{q/4} \right] + C_q t^{-q/2} \La 1, f^4 \Ra^{q/4} + C_q t^{-q/2} \La 1, f^2 \Ra^{q/2}
  \nonumber \\
   & \leq & C_q t^{-1} \La 1, f^4 \Ra^{q/4} + C_q t^{-q/2} \La 1, f^2 \Ra^{q/2}.
  \label{twoit}
\end{eqnarray}
One could repeat the trick more times here, but this estimate will already be sufficient for us. 
%
%
By normalising integrable $f,g$ we reach the following lemma.
\begin{Lemma} \label{decorr}
For $q \in 4 \N$, and conjugate $p$, there exists $c_6= c_6(q, \lambda_A, D_A,\mathcal{L}(\xi_0,\eta_0))$ so that
for integrable $f,g \geq 0$ 
\[
\left| \E \left[ \La \xi_t,f \Ra \La \eta_t,g \Ra \right] - \E \left[ \La \xi_t,f \Ra \right] \, 
\E\left[ \La \eta_t,g \Ra \right] \right| 
\leq c_6  \left( t^{-1/q} \La 1, f^4 \Ra^{1/4} + t^{-1/2} \La 1, f^2 \Ra^{1/2}\right) 
 \La g,1 \Ra \hat{\eta}_t^{1/p}.
\]
\end{Lemma}
%

\noindent
\textbf{One point estimates.}
Taking expectations in (\ref{mgpb}) and using translation invariance, we find
\begin{equation} \label{basicidea2}
d \hat{\eta}_t = - \lambda_B \E[ \xi_t(0) \eta_t(0)] dt
\end{equation} 
so that $t \to \hat{\eta}_t$ is decreasing. The test function formulation becomes
\begin{equation} \label{onepointB}
d \La \eta_t,\phi_t \Ra  =  \La \eta_t, \dot{\phi}_t + D_B \Delta \phi_t \Ra dt - \lambda_B 
\La \xi_t \eta_t, \phi_t \Ra dt + d \La M^B_t,\phi_t \Ra.
\end{equation}
Taking $\phi_{s}(x) = P^B_{t-s}f(x) := \sum_y p^B_{t-s}(y) f(x-y)$ for $s \in [0,t]$, we find
\begin{equation} \label{GreenB}
d \La \eta_s,P^B_{t-s}f \Ra = -  \lambda_B \La \xi_s \eta_s, P^B_{t-s}f \Ra ds + d \La M^B_s,P^B_{t-s}f \Ra.
\end{equation}
Taking expectations we find
\begin{eqnarray} \label{simpleboundB}
\left| \E[ \La \eta_t,f \Ra ] - \E[ \La \eta_{t-s}, P^B_s f \Ra] \right| & = & \lambda_B
\left| \int^t_{t-s} \E[ \La \xi_r \eta_r , P^B_{t-r}f \Ra ] dr \right|.
\end{eqnarray}
Without any correlation properties, to estimate the right hand side we will bound $\E[  \xi_r(0) \eta_r(0)]$ 
by developing a two point calculation. 

\noindent
\textbf{Two point estimates.}
Note that $d[\eta(x),\xi(y)]_t = 0$ since the processes never jump simultaneously. So calculus leads, for 
suitable test function $\phi_t(x,y)$, to
\begin{eqnarray}
d \La \xi_t * \eta_t,\phi_t \Ra & = & d \La \xi_t * \eta_t, \dot{\phi}_t + (D_A\Delta_x + D_B \Delta_y) \phi_t \Ra dt 
\nonumber \\
&& \;\; - \sum_x \sum_y \xi_t(x) \eta_t(y) \left( \lambda_B \xi_t(y) + \lambda_A (\xi_t(x)-1) \right) \phi_t(x,y) \, dt 
 + \mbox{MI}. \label{testfn3}
\end{eqnarray}
Define the test function
\begin{equation}
\psi^B_t(x,y) = \Pp[S^{A,x}_t=0, \,  S^{B,y}_t = 0,\, \tau >t]
\label{Btestfn}
\end{equation}
where $S^{A,x}$ is a rate $D_A$ simple random walks on $\Z^d$ starting at $x$,
$S^{B,y}$ is a rate $D_B$ simple random walks on $\Z^d$ starting at $y$, and
\[
\tau = \inf\{t: \int^t_0 I(S^{A,x}_s = S^{B,y}_s) ds > \mathcal{E}\}
\]
where $\mathcal{E}$ is an independent rate $\lambda_B$ exponential variable. 
Choosing $\phi_s=\psi^B_{t-s}$ for $s \in [0,t]$, we have $\phi_t(x,y) = I(x=y=0)$ and
\[
\dot{\phi}_s(x,y) + (D_A\Delta_x + D_B \Delta_y) \phi_s(x,y) = \lambda_B I(x=y) \phi_s(x,y), \quad \mbox{for $s \in [0,t]$.} 
\]
Using this in (\ref{testfn3}) we reach 
\begin{eqnarray}
d \La \xi_s * \eta_s,\phi_s \Ra &=&  - \lambda_B \sum_x \sum_{y: y \neq x} \xi_s(x) \xi_s(y) \eta_s(y) \phi_s(x,y) 
\, ds \nonumber \\
&& \;\; - \lambda_A \sum_x \sum_y  \xi_s(x)(\xi_s(x)-1) \eta_s(y) \phi_s(x,y) \, ds  + \mbox{MI}. 
\label{twopointB}
\end{eqnarray}
Note in particular, using $\psi_t(x,y) \leq p^A_t(x)p^B_t(y)$, that
\begin{equation} \label{AB2ptbound}
\E \left[ \xi_t(0) \eta_t(0) \right] \leq \E \left[ \La \xi_{t-s} * \eta_{t-s}, \psi_s \Ra \right]
\leq \E \left[ \La \xi_{t-s},p_s^A \Ra \La\eta_{t-s},p^B_s \Ra \right]. 
\end{equation}
Similarly, adjusting the test function, we can check 
\begin{equation} \label{AB2ptbound+}
\E \left[ \xi_t(x) \eta_t(y) \right] 
\leq \E \left[ \La \xi_{t-s},p_s^A(x-\cdot) \Ra \La\eta_{t-s},p^B_s(y-\cdot) \Ra \right].
\end{equation}
\begin{Lemma} \label{1pt}
For any $q \in 4\N $, with conjugate $p$, there exists $c_7 = c_7(\Law(\xi_0,\eta_0),\lambda_A,\lambda_B,D_A,q)$ 
so that whenever $0 \leq s \leq t/4$ and  $t \geq 2$, 
\begin{eqnarray*}
\left| \hat{\eta}_t - \hat{\eta}_{t-s} \right| \leq c_7 \, s \,
 \left( t^{-1} + s^{-(3d/8)} \right) 
\hat{\eta}_{t-2s}^{1/p}.
\end{eqnarray*}
\end{Lemma}

\textbf{Proof.} From (\ref{simpleboundB}) and (\ref{AB2ptbound}) we have
\begin{eqnarray}
\left| \hat{\eta}_t - \hat{\eta}_{t-s} \right| & = &   \lambda_B \int^t_{t-s} \E[ \xi_r(0) \eta_r(0)] dr \nonumber \\
& \leq & \lambda_B \int^t_{t-s} \E \left[ \La \xi_{t-2s},p_{r-t+2s}^A \Ra \La\eta_{t-2s},p^B_{r-t+2s} \Ra \right] dr.  \label{temp60}
\end{eqnarray}
Lemma \ref{decorr} with the choice $f = p_s^A$ and $g = p_s^B$ we find, for any $q \in 4 \N$ 
(and conjugate $p$), $s>0$ and $t \geq 1$,
\begin{eqnarray*}
\hspace{-.4in} \left| \E \left[ \La \xi_t,p_s^A \Ra  \La\eta_t,p_s^B \Ra \right] - \hat{\xi}_t \hat{\eta}_t \right| 
& \leq & c_6  \left( t^{-1/q} \La 1, (p^A_s)^4 \Ra^{1/4} + t^{-1/2} \La 1, (p^A_s)^2 \Ra^{1/2}\right) 
\hat{\eta}_t^{1/p} \\
& \leq & C \left( t^{-1/q} s^{-3d/8} + t^{-1/2} s^{-d/4} \right) 
\hat{\eta}_t^{1/p}
\end{eqnarray*}
using (\ref{pAbound}) in the final inequality. Using this in (\ref{temp60}) we get, using $r-t+2s \geq s$,  
\begin{eqnarray*}
\left| \hat{\eta}_t - \hat{\eta}_{t-s} \right| & \leq &  \lambda_B s \, 
\left( \hat{\eta}_{t-2s} \hat{\xi}_{t-2s} + C \left( (t-2s)^{-1/q} s^{-3d/8} + (t-2s)^{-1/2} s^{-d/4} \right) 
\hat{\eta}_{t-2s}^{1/p} \right) \\
& \leq & C s \left(  t^{-1} \hat{\eta}_{t-2s} + \left( t^{-1/q} s^{-3d/8} + t^{-1/2} s^{-d/4} \right) 
\hat{\eta}_{t-2s}^{1/p} \right) \quad \mbox{(using (\ref{crudeLUB}))}\\
& \leq & C \left(  s t^{-1} + s^{1-(3d/8)} \right) 
\hat{\eta}_{t-2s}^{1/p} 
\end{eqnarray*}
where we have thrown away $t^{-1/q}$, bounded $\hat{\eta}_t \leq C(p,\hat{\eta}_0) \hat{\eta}_{t}^{1/p}$ and
$ t^{-1/2} s^{1-(d/4)} \leq t^{-1} s + s^{1-(d/2)}$ for the final inequality. \qed

To improve the upper bound (\ref{AB2ptbound}) to an estimate we need to bound the terms on the right hand side
of (\ref{twopointB}), that is we need upper bounds on the three point terms
namely $\E[ \xi_t(x) \xi_t(y) \eta_t(y)]$ for $x \neq y$ and  $\E[ \xi_t(x)(\xi_t(x)-1) \eta_t(y)]$ for all $x,y$.

\noindent
\textbf{Three point estimates.} We are only looking for upper bounds. 
We will now show, analogously to (\ref{AB2ptbound+}), that when $0 \leq s \leq t$
\begin{eqnarray} 
\E \left[ \xi_t(x) \xi_t(y) \eta_t(y) \right] 
& \leq & \E \left[ \La \xi_{t-s},p_s^A(x-\cdot) \Ra \La \xi_{t-s},p_s^A(y-\cdot) \Ra \La\eta_{t-s},p^B_s(y-\cdot) \Ra \right] 
\label{3pt1} \\
\E \left[ \xi_t(x) (\xi_t(x)-1) \eta_t(y) \right] 
& \leq & \E \left[ \La \xi_{t-s},p_s^A(x-\cdot) \Ra^2 \La\eta_{t-s},p^B_s(y-\cdot) \Ra \right]. \label{3pt2}
\end{eqnarray}
Indeed using the test function $\phi_s(x,y) = \psi^{a,b}_{t-s}(x,y)$ from (\ref{psidefn}), we reach the decomposition
for $ d (\La \xi_s * \xi_s,\phi_s \Ra -  \La \xi_s,\tilde{\phi}_s \Ra)$ given in (\ref{39end}). We combine this with
the decomposition
\[
d \La \eta_s,p^B_{t-s}(\cdot- c) \Ra = -  \lambda_B \La \xi_s \eta_s, p^B_{t-s}(\cdot-c) \Ra ds + MI
\]
from (\ref{GreenB}). Note that $ \La \xi_s * \xi_s,\phi_s \Ra - \La \xi_s,\tilde{\phi}_s \Ra \geq 0$, 
$\La \eta_s,p^B_{t-s}(\cdot- c) \geq 0$ and 
these two processes have no simultaneous jumps. We see that the product
\[
\La \eta_s,p^B_{t-s}(\cdot- c) \Ra \left( \La \xi_s * \xi_s,\phi_s \Ra -  \La \xi_s,\tilde{\phi}_s \Ra \right) 
\]
is a positive supermartingale over $[0,t]$ since all the drift terms are negative. Comparing the expectations at times
$t$ and $t-s$ we reach 
\begin{eqnarray*}
&& \E[ \eta_t(c) \left( \xi_t(a) \xi_t(b) - \xi_t(a)I(a=b) \right) ] \\
& \leq &  \E\left[ \La \eta_{t-s}, p_s^B(\cdot -c) \Ra \left( \La \xi_{t-s} * \xi_{t-s}, \psi^{a,b}_s \Ra - \La \xi_{t-s}, \tilde{\psi}^{a,b}_s \Ra \right) \right] \\
& \leq &  \E \left[ \La \eta_{t-s}, p_s^B(\cdot -c) \Ra \left( \La \xi_{t-s} * \xi_{t-s}, p_s^A(\cdot-a) * p_s^A(\cdot -b) \Ra - 
\La \xi_{t-s}, p_s^A(\cdot-a) p_s^A(\cdot -b)  \Ra \right) \right] 
\end{eqnarray*}
where in the final inequality we used the bound $\psi^{a,b}_s (x,y) \leq p_s^A(x-a) p_s^A(y -b)$. Now taking
$c=b=y, \, a=x$ we reach (\ref{3pt1}), and taking $c=y, \, a=b=x$ we reach (\ref{3pt2}).

\begin{Lemma} \label{B3pt}
For $q \in 4 \N$, and conjugate $p$, there exist 
$c_8 = c_8(q,\Law(\xi_0,\eta_0),\lambda_A,D_A)$ so that for all $0<s\leq t/2$, $t \geq 1$ and $x,y$
\[
\E \left[ \xi_t(x) \xi_t(y) \eta_t(y) \right] \vee \E \left[ \xi_t(x) (\xi_t(x)-1) \eta_t(y) \right] 
\leq  c_8 \left( t^{-2} + s^{-3d/4} \right) \hat{\eta}_{t-s}^{1/p} 
\label{AB3ptbound2}
\]
\end{Lemma}

\noindent
\textbf{Proof.} We again apply H\"{o}lder's inequality to (\ref{3pt1},\ref{3pt2}) in order to bound them in terms of $\hat{\eta}$. 
We can write the right hand sides of (\ref{3pt1},\ref{3pt2}) as 
\begin{eqnarray}
&& \hspace{-.3in} \E \left[ \La \xi_{t-s},p_s^A(x-\cdot) \Ra \La \xi_{t-s},p_s^A(z-\cdot) \Ra \La\eta_{t-s},p^B_s(y-\cdot) \Ra \right] 
\nonumber \\
& = & \E \left[ \La \xi_{t-s} - \hat{\xi}_{t-s} ,p_s^A(x-\cdot) \Ra \La \xi_{t-s} - \hat{\xi}_{t-s},
p_s^A(z-\cdot) \Ra \La\eta_{t-s},p^B_s(y-\cdot) \Ra \right]  \nonumber \\
&& + \hat{\xi}_{t-s}  \E \left[ \La \xi_{t-s} - \hat{\xi}_{t-s} ,p_s^A(x-\cdot) \Ra  \La\eta_{t-s},p^B_s(y-\cdot) \Ra \right] \nonumber\\
&& + \hat{\xi}_{t-s}  \E \left[ \La \xi_{t-s} - \hat{\xi}_{t-s} ,p_s^A(z-\cdot) \Ra  \La\eta_{t-s},p^B_s(y-\cdot) \Ra \right] \nonumber\\
&& + \hat{\xi}_{t-s}^2 \E \left[ \La\eta_{t-s},p^B_s(y-\cdot) \Ra \right]. \label{temp31}
\end{eqnarray}
Arguing as in (\ref{StepA},\ref{StepB},\ref{twoit}) the first term of (\ref{temp31}) is bounded, when $q \in 4 \N$, by
\begin{eqnarray*}
&& \hspace{-.3in} \E [ \La \xi_{t-s}- \hat{\xi}_{t-s},p_s^A \Ra^{2q} ]^{1/q} 
\, \E [ \La  \tilde{\eta}_{t-s},p^B_s \Ra ]^{1/p} \\
& \leq & C  \left((t-s)^{-1} \La 1,(p_s^A)^4 \Ra^{q/2} + (t-s)^{-q} \La 1,  (p_s^A)^2 \Ra^q \right)^{1/q}
 \hat{\eta}_{t-s}^{1/p} \\
& \leq & C\left(t^{-1/q} s^{-3d/4} + t^{-1} s^{-d/2} \right) \hat{\eta}_{t-s}^{1/p}
\end{eqnarray*}
when $0 < s \leq t/2$ and $t \geq 1$; here $C = C(q,(\xi_0,\eta_0), D_A,\lambda_A)$. 
Similarly, the second and third term of (\ref{temp31}) are bounded by 
$C t^{-1} \left( t^{-1/q} s^{-3d/8} + t^{-1/2} s^{-d/4} \right) \hat{\eta}_{t-s}^{1/p}$.
The final term of (\ref{temp31}) equals $\hat{\xi}_{t-s}^2 \hat{\eta}_{t-s} \leq C t^{-2} \hat{\eta}_{t-s}$. 
Collecting these bounds we reach, after grouping as in Lemma \ref{1pt}, the stated estimate. 
\qed

These bounds can be used in our final estimate, which continues from the expansion in (\ref{twopointB}) using the test function
$\psi^B$ defined in (\ref{Btestfn}).
\begin{Lemma} \label{B2pt}
For $q \in 4 \N$, and conjugate $p$, there
there exists 
$c_9 = c_9(\Law(\xi_0,\eta_0),\lambda_A, \lambda_B,D_A,D_B,q)$ so that  
for $0<s<t/4$ and $t \geq 2$ 
\[
\left| \E \left [\xi_t(0)\eta_t(0) \right] - \E \left[ \La \xi_{t-s} * \eta_{t-s},\psi^B_s \Ra \right] \right| \leq 
c_9 \, s \,\left( t^{-2} + s^{-3d/4} \right)  \hat{\eta}_{t-2s}^{1/p}. 
\]
\end{Lemma}
\noindent
\textbf{Proof.} Integrate (\ref{twopointB}) over the interval $[t-s,t]$. The estimate will follow if we bound the expectation of the two
drift terms
\begin{eqnarray*}
&& \hspace{-.4in} \lambda_B \int^t_{t-s} \sum_x \sum_{y: y \neq x} \xi_{r}(x) \xi_{r}(y) \eta_{r}(y) \psi^B_{t-r}(x,y) dr \\
&& + \lambda_A \int^t_{t-s}\sum_x \sum_{y} \xi_{r}(x)(\xi_{r}(x)-1) \eta_{r}(y)  \psi^B_{t-r}(x,y) dr.
\end{eqnarray*}
We use $\psi^B_s(x,y) \leq p^A_s(x) p^B_s(y)$ and then the bounds from Lemma \ref{B3pt} to estimate this in terms of
$\hat{\eta}_{t-2s}$; the requirement that $s \leq t/4$ and $t \geq 2$ ensure this lemma is applicable. \qed

Finally we need an estimate analogous to (\ref{KdBl12}) for the non-coalescence
probability $\psi^B_t(x,y)$. Namely, when $d \geq 3$,
for some $c_{10} = c_{10}(D_A, D_A,d)$
\begin{equation} \label{KdBl12.2}
\sum_{x,y \in\Z^d} \left| \psi^B_t(x,y) - p_B p^A_t(x) p^B_t(y) \right| \leq c_{10} \, t^{-\delta} \quad 
\mbox{for all $t \geq 0$}.
\end{equation}
The argument from Lemma 12 
in \cite{kesten1998asymptotic} applies for two particles with different jumps rates $D_A$ and $D_B$,
and again the estimate holds for any $0<\delta < \frac12 \frac{(d-2)}{(d-1)}$.

\subsection{$B$-particle modified rate equation} \label{s4.4}
We repeat the main approximations from the sketch in the introduction:
\begin{eqnarray}
\frac{d \hat{\eta}_t}{dt} & = &  - \lambda_B \E[\xi_t(0)\eta_t(0)] \nonumber \\
& \approx & - \lambda_B \E[\La \xi_{t-s} * \eta_{t-s}, \psi^B_s \Ra]  \label{app5} \\
& \approx & - \lambda_B p_B  \E[ \La \xi_{t-s},p^A_s \Ra \La \eta_{t-s},p^B_s \Ra  ]  \label{app6} \\
& \approx & - \lambda_B p_B  \E[ \La \xi_{t-s},p^A_s \Ra] E[ \La \eta_{t-s},p^B_s \Ra  ]  
\label{app7} \\
& = & - \lambda_B p_B  \hat{\xi}_{t-s} \hat{\eta}_{t-s} \nonumber \\
& \approx & - \lambda_B p_B  \hat{\xi}_{t} \hat{\eta}_{t} \label{app8} \\
& \approx & -  \frac{\lambda_B p_B}{\lambda_A p_A } \frac{1}{t} \hat{\eta}_{t}. \label{app9} 
\end{eqnarray}
We will check below that the lemmata from Section \ref{s4.3} quantify these approximations and yield
\begin{equation} \label{ereB}
\frac{d \hat{\eta}_t}{dt}
=  - \theta t^{-1} \hat{\eta}_t + \mathcal{E}_t
\end{equation}
where $\theta = \lambda_B p_B / \lambda_A p_A$
and, for some $\kappa >0, \, t_1 \geq 2$ and any $p>1$, 
\begin{equation} \label{Bparticleerror}
\left| \mathcal{E}_t \right| \leq C t^{-(1+\kappa)}  \hat{\eta}_{t-2s}^{1/p}   \quad \mbox{when $0<s \leq t/4$ and $t \geq t_1$}
\end{equation}
with $C = C(\Law(\xi_0,\eta_0),\lambda_A, \lambda_B,D_A,D_B,p)$.
To analyse (\ref{ereB}) let
\[
t_* = \inf\{t \geq 0: \hat{\eta}_t \geq K t^{-\theta}\}.
\]
By choosing $K$ large we can ensure $t_* \geq t_1$. 
We will argue that $t_* = \infty$ if we choose $K$ suitably.
Note that for $t_1 \leq  t \leq t_*$
\begin{eqnarray*}
 t^{\theta} \mathcal{E}_t 
& \leq &  C t^{\theta}  t^{-(1+\kappa)} \hat{\eta}_{t-2s}^{1/p}  \\
& \leq & C K^{1/p} t^{\theta}  t^{-(1+\kappa)} (t-2s)^{-\theta/p} \\
& \leq & C K^{1/p} t^{\theta -(1+\kappa)-(\theta/p)} 
\end{eqnarray*}
using $s \leq t/4$ in the final inequality (and letting $C$ vary line to line).
Choose $p>1$ so that $\theta- (1+\kappa) -(\theta/p)  < - 1$.
Suppose $t_1 \leq t_* < \infty$. Then, integrating  $d (t^{\theta} \hat{\eta}_t)/dt = t^{\theta} \mathcal{E}_t $ over $[t_1,t_*]$, we have
\begin{eqnarray*}
K & = & t_*^{\theta} \, \hat{\eta}_{t_*} \\
& = & t_1^{\theta} \hat{\eta}_{t_1}  + \int^{t_*}_{t_1} s^{\theta} \mathcal{E}_s ds \\
& \leq & t_1^{\theta} \hat{\eta}_{t_1}   + C K^{1/p}
 \int^{t_*}_{t_1} s^{\theta- (1+\kappa) -(\theta/p)}  ds. 
\end{eqnarray*}
This is a contradiction if we choose $K=K(\theta,p,t_1,\hat{\eta}_{t_1})$ large enough. 
We conclude $t_* = \infty$ for suitable $K$ and then we
integrate up to find
\[
t^{\theta} \hat{\eta}_t = {t_1}^{\theta} \hat{\eta}_{t_1}  + \int^t_{t_1} s^{\theta} \mathcal{E}_s ds  
= t_1^{\theta} \hat{\eta}_{t_1} + \int^{\infty}_{t_1} s^{\theta} \mathcal{E}_s ds  + O(t^{\theta - \kappa - (\theta/p)})
\]
establishing the desired asymptotics in Theorem \ref{mainresult} with $\kappa_2 = \kappa -\theta + (\theta/p)$. 
Note that by taking $p$ large we can take $\kappa_2$ as close to $\kappa$ as desired. 

It remains to verify the error bound stated in (\ref{Bparticleerror}). We will choose $s = t^{\alpha}$, for some
$\alpha \in (0,1)$ shortly and then we choose $t_1$ so that $t \geq t_1$ implies that $t \geq 2$ and $s \leq t/4$. 
Throughout $C$ may depend on $p,d,\lambda_A,\lambda_B, D_A,D_B, \mathcal{L}(\xi_0,\eta_0)$. 

The error in  (\ref{app5}) is controlled by Lemma \ref{B2pt} and is bounded by 
$C s \,\left( t^{-2} + s^{-3d/4} \right)  \hat{\eta}_{t-2s}^{1/p}$. 
The error in  (\ref{app6}) is bounded, using  (\ref{KdBl12.2}) by 
\[
\lambda_B \sum_{x,y} \E[ \xi_{t-s}(x) \eta_{t-s}(y)] \left| \psi^B_s(x,y) - p_B p^A_s(x) p^B_s(y) \right|
\leq C s^{-\delta} \sup_{x,y} \E[ \xi_{t-s}(x) \eta_{t-s}(y)].
\]
Using (\ref{AB2ptbound+}) and Lemma \ref{decorr} once more we have
\begin{eqnarray*}
\E \left[ \xi_{t-s}(x) \eta_{t-s} (y) \right] & \leq & \E \left[ \La \xi_{t-2s},p_s^A(x-\cdot) \Ra \La \eta_{t-2s},p^B_s(y-\cdot) \Ra \right] \\
& \leq & \hat{\xi}_{t-2s} \hat{\eta}_{t-2s} + C \left( t^{-1/q} \La 1, (p_s^A)^4 \Ra^{1/4} + t^{-1/2} \La 1, (p_s^B)^2 \Ra \right) 
\hat{\eta}^{1/p}_{t-2s} \\
& \leq & Ct^{-1}  \hat{\eta}_{t-2s} + C \left( s^{-3d/8} + t^{-1/2} s^{-d/4} \right) \hat{\eta}^{1/p}_{t-2s} \\
& \leq  & C \left( t^{-1} + s^{-3d/8}\right) \hat{\eta}^{1/p}_{t-2s}.
\end{eqnarray*}
Thus the error in (\ref{app6}) is bounded by $C s^{-\delta} \left( t^{-1} + s^{-3d/8}\right)  \hat{\eta}^{1/p}_{t-2s}$. 
The error in  (\ref{app7}) is the key covariance estimate in Lemma
\ref{decorr} and, as above, is bounded by 
$ C \left(s^{-3d/8} + t^{-1/2} s^{-d/4} \right) \hat{\eta}_{t-s}^{1/p}$. 
For the error in (\ref{app8}) we combine Lemma \ref{A1pt} and Lemma \ref{1pt} to obtain the bound 
$C s t^{-1} \left( t^{-1}  + s^{-(3d/8)} \right) \hat{\eta}_{t-2s}^{1/p} $. 
The final error in  (\ref{app9}) arises from the A-particle asymptotics in Theorem \ref{mainresult} and yields a term
$C\hat{\eta}_t t^{-1-\kappa_1}$.
Combining all the errors, using that $t \to \hat{\eta}_t$ is decreasing as always, and omitting terms that
are dominated by other terms, we reach, for $0<s \leq t/4$ and $t \geq 2$, 
\[
\left| \mathcal{E}_t \right| \leq C \left(st^{-2} + s^{-3d/8} + t^{-1/2} s^{-d/4} + s^{-\delta} t^{-1} + t^{-1-\kappa_1} \right) \hat{\eta}_{t-2s}^{1/p}   \quad \mbox{when $0<s \leq t/4$ and $t \geq t_1$}
\]
Now choosing $s=t^{\alpha}$ for any $\alpha \in (8/3d,1)$ produces the required estimate in  (\ref{Bparticleerror}).

\noindent
\textbf{Remark 4.4.1} The optimal choice of $\alpha$ is a bit fiddly and $d$ dependent. In $d=3$ by choosing $\alpha = \frac{16}{17}$ one achieves $\kappa_2 = \frac{1}{17}$; as $d \to \infty$ the error term $t^{-1-\kappa_1}$ dominates, so that 
$\kappa_2 = \kappa_1$ for large $d$.

\noindent
\textbf{Remark 4.4.2} We make brief comments about the models where
either $\lambda_A$ or $\lambda_B$, or both, are infinite. One method to establish 
the correct decay rates would be to pass to the limits $\lambda_A \to \infty$ or
$\lambda_B \to \infty$ in the modified rate equations. The finite rate models, which we might denote by
 $\xi^{\lambda_A}$ and $\eta^{\lambda_A,\lambda_B}$, converge (at leats at fixed $t$) to their
 infinite rate counterparts $\xi^{\infty}$ and $\eta^{\lambda_A, \infty}$ or $\eta^{\infty, \infty}$. 
 Passing to the limit $\lambda_A \to \infty$ in (\ref{Aere}), or the limit $\lambda_B \to \infty$
 in (\ref{ereB}) we expect to reach
 \[
 \frac{d \hat{\xi}^{\infty}_t}{dt} =  - \gamma D_A  (\hat{\xi}^{\infty}_{t})^2 + \mathcal{E}_t, 
 \quad \mbox{or} \quad \frac{d \hat{\eta}^{\infty,\infty}_t}{dt}
=  - \theta t^{-1} \hat{\eta}^{\infty,\infty}_t + \mathcal{E}_t
 \]
 with the appropriate limiting value of $\theta$. The catch is that we need to bound
the errors $\mathcal{E}_t$ from (\ref{Aere},\ref{ereB}) uniformly over large $\lambda_A$ or $\lambda_B$.
 This is immediate for terms such as (\ref{app2},\ref{app3},\ref{app4}) which are bounded using the 
 product $\lambda_A p_A$, which is bounded in $\lambda_A$, but looks less easy for
 the approximation (\ref{app1}). 
%
%
 We believe it is easier to restart the entire argument using the equations for the infinite rate models;
the starting point in the case where both rates are infinite are the exact formulae
\[
\frac{d \hat{\xi}_t}{dt}  =   - \E[\xi_t(0) \xi_t(e)], \qquad 
\frac{d \hat{\eta}_t}{dt} =   - \E[\xi_t(0)\eta_t(e)]
\]
where $e$ is a neighbouring site to the origin. Indeed an $A$ particle modified
rate equation for the case of instantaneous coalescence was established in $d=2$
(where extra logarithms emerge) in \cite{lukins2018multi} starting this way. 
\section{Appendix} 
We show here how the estimates on $A$ particles in Section \ref{s4.2} can be used to derive the upper bound in (\ref{crudeLUB}). We repeat the steps (\ref{app1},\ref{app2}), estimating errors only from above, as follows: for $0<s<t$
\begin{eqnarray}
\frac{d \hat{\xi}_t}{dt} & = &  - \lambda_A \E[\xi_t(0)(\xi_t(0)-1)] \nonumber \\
& \leq & - \lambda_A \E[\La \xi_{t-s} * \xi_{t-s}, \psi^{0,0}_s \Ra]  + \mathcal{E}_1 \label{app1+} \\
& \leq & - \lambda_A p_A  \E[ \La \xi_{t-s},p^A_s \Ra^2] + \mathcal{E}_1 +  \mathcal{E}_2 
\label{app2+} \\ 
& \leq & - \lambda_A p_A  (\E[ \La \xi_{t-s},p^A_s \Ra])^2  + \mathcal{E}_1 +  \mathcal{E}_2  
 \nonumber \\
&   =  &- \lambda_A p_A  \hat{\xi}_{t-s}^2 + \mathcal{E}_1 +  \mathcal{E}_2.  \nonumber 
\end{eqnarray}
The error $\mathcal{E}_2$ can be estimated using (\ref{KdBl12}), and $\psi^{0,0}_t(x,y) \leq  p^A_t(x) p^A_t(y)$,  by
\begin{eqnarray}
\mathcal{E}_2 & \leq & \lambda_A  \sum_{x,y}  \E[  \xi_{t-s}(x) \xi_{t-s}(y)] \left| p_A p^A_s(x) p^A_s(y)   - \psi^{0,0}_s(x,y) \right| \nonumber \\
& \leq &  \lambda_A  \hat{\xi}_{t-s}^2 \sum_{x \neq y}  \left| p_A p^A_s(x) p^A_s(y)   - \psi^{0,0}_s(x,y) \right| + \lambda_A \sum_x \E[ \xi_{t-s}^2(x)] (p^A_s(x))^2 \nonumber  \\
& \leq & C(D_A,\lambda_A) \left( \hat{\xi}_{t-s}^2 s^{-\delta} + (\hat{\xi}_{t-s} + 2 \hat{\xi}^2_{t-s}) s^{-d/2} \right).\label{E2}
\end{eqnarray}
Here we have used the simple bound (\ref{pAbound}), we have bounded the second 
moment using $\E[ \xi_{t-s}^2(x)] = \E[ \xi_{t-s}(x)(\xi_{t-s}(x)-1)] + \E[ \xi_{t-s}^2(x)]$, and used 
negative correlation (\ref{natwo}). To bound $\mathcal{E}_1$ we revisit the three terms in (\ref{temp107}).
Arguing as in Lemma \ref{A2pt} (see (\ref{pAbound}) and (\ref{41.5})),
\begin{eqnarray*}
\int^t_{t-s}  \sum_x \sum_{y: y \neq x}  \E \left [ \xi_r(x) \xi_r(y) (\xi_r(x)+\xi_r(y)-2) \right]  \phi_r(x,y)  dr & \leq & 4 s \,  \hat{\xi}_{t-s}^3; 
 \end{eqnarray*}
and
\begin{eqnarray*}
\E \left [\La \xi_{t-s}, \tilde{\phi}_{t-s} \Ra \right] & \leq & C(D_A) \hat{\xi}_{t-s} s^{-d/2}
\end{eqnarray*}
 and
 \begin{eqnarray*}
\int^t_{t-s} \E \left [\La \xi_r(\xi_r-1)(\xi_r-2),\tilde{\phi}_r \Ra \right]
  & \leq & C(D_A) \hat{\xi}_{t-s}^3 \int^s_0 (r^{-d/2} \wedge 1).
\end{eqnarray*}
Together these imply using  (\ref{temp107}) that for all $0<s <t$
\begin{equation} \label{E1}
\mathcal{E}_1 \leq C(\lambda_A,D_A) \left( (1+s)  \hat{\xi}_{t-s}^3 +  s^{-d/2} \hat{\xi}_{t-s} \right).
\end{equation}
Using (\ref{E1}) and (\ref{E2}) we reach, for any $0<s<t$ and $c_0 = c_0(\lambda_A,D_A)$, 
\begin{equation} \label{Meqn}
\frac{d \hat{\xi}_t}{dt} \leq - \hat{\xi}_{t-s}^2 ( \lambda_A p_A  - c_0s^{-\delta} - c_0s^{-d/2}) + c_0 
(1+s)  \hat{\xi}_{t-s}^3 +  c_0 s^{-d/2} \hat{\xi}_{t-s}.
\end{equation}
The aim is to choose $s$ carefully to show that the linear and cubic terms are dominated by the 
quadratic term. 
We will choose $s = s(t) \in [0,t]$ to be the solution to the implicit equation
\begin{equation} \label{fixedpoint}
s = \alpha \hat{\xi}^{-1}_{t-s}, 
\end{equation}
where we will choose $\alpha= \alpha(\lambda_A,D_A) >0$ shortly. Note that, for a fixed $t>0$, 
the function $s \to \alpha \hat{\xi}^{-1}_{t-s}$ is decreasing, so that there is a unique solution
$s(t)$ as soon as $t \geq t_0 := \alpha \hat{\xi}^{-1}_{0}$. We note that $t \to s(t) $ is increasing and we
claim that
\begin{equation} \label{claims}
\mbox{(i) $s(t) \to \infty$ as $t \to \infty$ \hspace{.2in} and \hspace{.2in}  (ii) $t-s(t) \to \infty$ as $t \to \infty$.}
\end{equation}
To see (i) we use $\hat{\xi}_t \downarrow 0$ as $ t \to \infty$; hence for any $M$ we may find
$t_M$ so that $\alpha/\hat{\xi}_{t_M} = M$ and then $s(t_M+M) = M$. To see (ii) we use the lower bound 
$\hat{\xi_t} \geq \left( \hat{\xi}_0^{-1} + \lambda_A t \right)^{-1}$ from (\ref{LB}) so that
$s(t) = \alpha \hat{\xi}^{-1}_{t-s(t)} \leq \alpha (\hat{\xi}_0^{-1} + \lambda_A (t-s(t))$.

Now we take $t_1 \geq t_0$ so that for $t \geq t_1$ we have $s(t) \geq 1$ and 
\[
 \lambda_A p_A  - c_0s^{-\delta}(t) - c_0s^{-d/2}(t) \geq \frac12 \lambda_A p_A.
 \]
 Using $s=s(t) $ in (\ref{Meqn}) we find
\begin{eqnarray*} \frac{d \hat{\xi}_t}{dt} & \leq & - \hat{\xi}_{t-s}^2  \frac{\lambda_A p_A}{2}  + 2 c_0 s  \hat{\xi}_{t-s}^3 +  c_0 \hat{\xi}_{t-s} s^{-d/2} \\
& = & - \hat{\xi}_{t-s}^2  \left( \frac{\lambda_A p_A}{2} - 2 c_0 \alpha \right) +  c_0 \alpha^{-d/2}
\hat{\xi}^{1+(d/2)}_{t-s}.
\end{eqnarray*}
Now we choose $\alpha $ so that $ \frac{\lambda_A p_A}{2} - 2 c_0 \alpha = \frac{\lambda_A p_A}{3}$. 
Since $\hat{\xi}_t \to 0$ and $t-s(t) \to \infty$ as $t \to \infty$  we may choose $t_2\geq t_1$ so that 
\[
c_0 \alpha^{-d/2}
\hat{\xi}^{1+(d/2)}_{t-s(t)} \leq \frac{\lambda_A p_A}{4} \hat{\xi}^{2}_{t-s(t)} \quad \mbox{for $t \geq t_2$.}
\]
We have reached
\[
\frac{d \hat{\xi}_t}{dt} \leq  - \frac{\lambda_A p_A}{12} \hat{\xi}_{t-s}^2   \leq - \frac{\lambda_A p_A}{12} \hat{\xi}_{t}^2 \quad \mbox{for $t \geq t_2$.}
\]
which implies the desired upper bound. 

\bibliographystyle{abbrv}
\bibliography{AB_bib}

\begin{thebibliography}{10}

\bibitem{lukins2017thesis}


\bibitem{arratia1981limiting}
R.~Arratia.
\newblock Limiting point processes for rescalings of coalescing and
  annihilating random walks on zd.
\newblock {\em The Annals of Probability}, pages 909--936, 1981.

\bibitem{bramson1980asymptotics}
M.~Bramson and D.~Griffeath.
\newblock Asymptotics for interacting particle systems on $z^d$.
\newblock {\em Zeitschrift f{\"u}r Wahrscheinlichkeitstheorie und verwandte
  Gebiete}, 53(2):183--196, 1980.

\bibitem{chow2003probability}
Y.~S. Chow and H.~Teicher.
\newblock {\em Probability theory: independence, interchangeability,
  martingales}.
\newblock Springer Science \& Business Media, 2003.

\bibitem{dynkin1965markov}
E.~Dynkin.
\newblock {\em Markov processes}.
\newblock Springer, 1965.

\bibitem{garrod2018examples}
B.~Garrod, M.~Poplavskyi, R.~P. Tribe, and O.~V. Zaboronski.
\newblock Examples of interacting particle systems on z as pfaffian point
  processes: annihilating and coalescing random walks.
\newblock In {\em Annales Henri Poincar{\'e}}, volume~19, pages 3635--3662.
  Springer, 2018.

\bibitem{graham1996weak}
C.~Graham, T.~G. Kurtz, S.~M{\'e}l{\'e}ard, P.~E. Protter, M.~Pulvirenti,
  D.~Talay, T.~G. Kurtz, and P.~E. Protter.
\newblock Weak convergence of stochastic integrals and differential equations
  ii: Infinite dimensional case.
\newblock {\em Probabilistic Models for Nonlinear Partial Differential
  Equations: Lectures given at the 1st Session of the Centro Internazionale
  Matematico Estivo (CIME) held in Montecatini Terme, Italy, May 22--30, 1995},
  pages 197--285, 1996.

\bibitem{howard}
M.~Howard.
\newblock Fluctuation kinetics in a multispecies reaction-diffusion system.
\newblock {\em Journal of Physics A: Mathematical and General}, 29(13):3437,
  1996.

\bibitem{kesten1998asymptotic}
H.~Kesten and J.~van~den Berg.
\newblock Asymptotic density in a coalescing random walk model.
\newblock {\em CWI. Probability, Networks and Algorithms [PNA]}, (R 9815),
  1998.

\bibitem{redner}
P.~Krapivsky, E.~Ben-Naim, and S.~Redner.
\newblock Kinetics of heterogeneous single-species annihilation.
\newblock {\em Physical Review E}, 50(4):2474, 1994.

\bibitem{key1957219m}
T.~M. Liggett.
\newblock Negative correlations and particle systems.
\newblock {\em Markov Process. Related Fields}, 8(4):547--564, 2002.
\newblock MR:1957219. Zbl:1021.60084.

\bibitem{lukins2018multi}
J.~Lukins, R.~Tribe, and O.~Zaboronski.
\newblock Multi-point correlations for two-dimensional coalescing or
  annihilating random walks.
\newblock {\em Journal of Applied Probability}, 55(4):1158--1185, 2018.

\bibitem{monthus}
C.~Monthus.
\newblock Exponents appearing in heterogeneous reaction-diffusion models in one
  dimension.
\newblock {\em Physical Review E}, 54(5):4844, 1996.

\bibitem{newman1984asymptotic}
C.~M. Newman.
\newblock Asymptotic independence and limit theorems for positively and
  negatively dependent random variables.
\newblock {\em Lecture Notes-Monograph Series}, pages 127--140, 1984.

\bibitem{pemantle2000towards}
R.~Pemantle.
\newblock Towards a theory of negative dependence.
\newblock {\em Journal of Mathematical Physics}, 41(3):1371--1390, 2000.

\bibitem{rajesh}
R.~Rajesh and O.~Zaboronski.
\newblock Survival probability of a diffusing test particle in a system of
  coagulating and annihilating random walkers.
\newblock {\em Physical Review E—Statistical, Nonlinear, and Soft Matter
  Physics}, 70(3):036111, 2004.

\bibitem{reimer2000proof}
D.~Reimer.
\newblock Proof of the van den berg--kesten conjecture.
\newblock {\em Combinatorics, Probability and Computing}, 9(1):27--32, 2000.

\bibitem{tz3}
R.~Tribe and O.~Zaboronski.
\newblock Pfaffian formulae for one dimensional coalescing and annihilating
  systems.
\newblock {\em Electronic Journal of Probability}, 16:2080--2103, 2011.

\bibitem{van2002randomly}
J.~van~den Berg and H.~Kesten.
\newblock {\em In and out of equilibrium, vol. 51 of Progress in Probability;
  Randomly coalescing random walk in dimension $d \geq 3$}.
\newblock Birkhauser, Boston, 2002.

\end{thebibliography}

\end{document}